\documentclass[12pt]{amsart}
\usepackage{amsmath,amssymb,amsthm,mathtools,mathrsfs}
\usepackage[T1]{fontenc}
\usepackage[colorlinks,linkcolor=red,citecolor=blue,urlcolor=blue]{hyperref}
\newtheorem{thm}{Theorem}[section]
\newtheorem{prop}[thm]{Proposition}
\newtheorem{lem}[thm]{Lemma}
\newtheorem{cor}[thm]{Corollary}
\theoremstyle{definition}

\theoremstyle{remark}
\newtheorem{rem}[thm]{Remark}
\numberwithin{equation}{section}
\newcommand{\abs}[1]{\left\lvert#1\right\rvert}
\newcommand{\norm}[1]{\left\|#1\right\|}
\renewcommand{\geq}{\geqslant}
\renewcommand{\leq}{\leqslant}
\newcommand{\e}{\varepsilon}
\newcommand{\R}{\mathbb{R}}
\newcommand{\C}{\mathbb{C}}
\newcommand{\Q}{\mathbb{Q}}

\newcommand{\mP}{\mathbb{P}}
\newcommand{\cM}{\mathcal{M}}
\newcommand{\cO}{\mathcal{O}}
\newcommand{\cE}{\mathcal{E}}
\newcommand{\cH}{\mathcal{H}}
\newcommand{\cR}{\mathcal{R}}
\newcommand{\cX}{\mathcal{X}}
\newcommand{\cY}{\mathcal{Y}}
\newcommand{\cL}{\mathcal{L}}
\newcommand{\cA}{\mathcal{A}}
\newcommand{\NA}{\mathrm{NA}}
\newcommand{\rad}{\mathrm{rad}}
\newcommand{\red}{\mathrm{red}}
\newcommand{\ddbar}{\frac{\sqrt{-1}}{2\pi}\partial\bar\partial}

\DeclareMathOperator{\MA}{MA}
\DeclareMathOperator{\Cal}{Cal}
\DeclareMathOperator{\DF}{DF}

\DeclareMathOperator{\Hom}{Hom}
\DeclareMathOperator{\Ric}{Ric}
\DeclareMathOperator{\tr}{tr}
\DeclareMathOperator{\vol}{vol}
\DeclareMathOperator{\PSH}{PSH}
\DeclareMathOperator{\ord}{ord}
\DeclareMathOperator{\supp}{supp}
\DeclareMathOperator{\RePart}{Re}
\DeclareMathOperator{\diag}{diag}

\DeclareMathOperator*{\argmin}{argmin}
\title[Calabi functional and non-Archimedean K-energy]{Lower bound of the Calabi functional\protect\\and non-Archimedean K-energy}
\author{Tomoyuki HISAMOTO}
\date{October 11, 2026}


\begin{document}
\begin{abstract}
We prove Donaldson's conjecture which describes the infimum
of the Calabi functional by the normalized
non-Archimedean K-energy.
The proof makes use of 
proximal equation for the K-energy. 
To obtain approximate smooth solutions, a continuity path of prescribed scalar curvature equation is coupled with a version of Chen's $\e$-subgeodesic. 
It does not require a global smooth solution
to the Calabi flow.

\end{abstract}
\maketitle
\setcounter{tocdepth}{1}
\tableofcontents

\section{Introduction}

The functional in the title was introduced by Calabi in his seeking a notion of standard metric form a variational viewpoint \cite{Cal82}.  
Donaldson later drew an infinitesimal moment map picture in \cite{Don97} and subsequently gave a 
lower bound of the Calabi functional in \cite{Don05}. 
The lower bound is obtained by 
his generalization of Futaki invariants \cite{Fut83,Fut83cscK}
and it was further predicted that this lower bound is tight. 

To make the statement more precisely, let us prepare terminology.  
Throughout the paper, $X$ is a compact K\"ahler manifold of dimension $n\geq1$, and $L$ is an ample line bundle. For a K\"ahler metric $\omega\in c_1(L)$, we set the volume as $V:=\int_X\omega^n=L^n$ and the Ricci curvature as 
\begin{equation*}
\Ric(\omega):= -\ddbar \log\omega^n 
=-dd^c\log\omega^n. 
\end{equation*}
The scalar curvature $S(\omega)=\tr_\omega\Ric(\omega)$ is simply obtained by taking trace. 
The average is computed as $\bar S=-nV^{-1}(K_XL^{n-1})$. 
The Calabi functional 
\begin{equation}\label{eq:calabi}
 \Cal(\omega):=\left[\frac{1}{V}\int_X(S(\omega)-\bar S)^2\omega^n\right]^{\frac{1}{2}}.
\end{equation}
measures how a given metic is far from the constant scalar curvature K\"ahler (shortly, cscK) condition. 
If one fixes the reference $\omega_0$, by 
$dd^c$-lemma, up to addition of constants, 
collection of $\omega$ is naturally identified with   
the space of smooth K\"ahler potential 
\begin{equation}\label{eq:cH}
 \cH:=\bigg\{\varphi\in C^\infty(X,\R): \omega_\varphi:=\omega_0+dd^c\varphi>0\bigg\}.
\end{equation}
As a consequence of this description, each tangent space of $\cH$ is naturally 
identified with the space of 
$C^\infty$ functions. 

The functional \eqref{eq:calabi} is closely related to K-energy functional defined by its differential: 
\begin{equation}\label{eq:mabuchi}
(d_{\varphi} M) f = - \frac{1}{V}\int_X f(S(\omega_{\varphi})-\bar{S}) \omega_{\varphi}^n. 
\end{equation}
K-energy was introduced by Mabuchi \cite{Mab86}, as an integration of 
Futaki's obstruction to the K\"ahler--Einstein metric.  
Based on works of Futaki and Mabuchi, 
Ding--Tian \cite{DT92,Tian97} generalize Futaki's obstruction 
to infinitesimal slope of the K-energy  
along a certain degeneration of a Fano manifold. Donaldson \cite{Don02} introduced 
the notion of K-stability, giving 
a purely algebraic definition 
of this generalized Futaki invariant 
for an arbitrary polarized manifold. 
We follow the 
non-Archimedean treatment of \cite{BHJ17} and formulate Donaldson's conjecture using the non-Archimedean K-energy $M^{\NA}(\cX,\cL)$ \cite[Definition 7.13]{BHJ17} of a test configuration $(\cX,\cL)$. 
A polarized manifold $(X, L)$ is called 
{\em K-semistable} if $M^{\NA}(\cX, \cL) \geq 0 $ holds for all test configurations $(\cX,\cL)$.  
Otherwise it is called {\em K-unstable}. 
We also need the $L^2$-norm 
\begin{equation}\label{eq:tc-norm}
 \norm{(\cX,\cL)}_2
 :=\left[\int_\R(\lambda-\bar\lambda)^2d\mu_{(\cX,\cL)}(\lambda) \right]^\frac{1}{2}
\end{equation}
of a test configuration. 
The right-hand side is square root of the variance of the probability Duistermaat--Heckman measure. A precise definition will be given in Section~\ref{sec:normalization}.
Let us firstly state our main result. 

\begin{thm}\label{thm:main}
For every K-unstable smooth complex projective polarized manifold $(X,L)$,
\begin{equation}\label{eq:main}
 \inf_{\varphi\in\cH} \Cal(\omega_{\varphi})
 =\sup_{(\cX, \cL)}
 \frac{-M^{\NA}(\cX,\cL)}{\norm{(\cX,\cL)}_2}.
\end{equation}
Here the supremum is taken over normal ample algebraic test configurations with non-zero norms. 
\end{thm}

Donaldson proved the one-sided inequality $\geq$ and predicted that 
the equality holds. 
It also takes an analogue of the Yang-Mills functional 
over a Riemann surface \cite{AB83}. 
One can further expects that the supremum is actually the maximum, 
however, we do not enter this harder problem in this article. 

Let us look back earlier research around Theorem \ref{thm:main}.  
The equality \eqref{eq:main} is modelled on so-called 
{\em moment--weight inequality} in classical 
Kempf-Ness theory, which relate the square norm of the moment map to  GIT-instability. 
Donaldson's proof actually starts from  these  finite-dimensional moment maps and uses the asymptotic expansion of the Bergman kernel in order to link them with the inequality for the infinite-dimensional space. In the toric setting, Sz\'ekelyhidi constructed an optimal destabilizing convex function and proved the corresponding sharp Calabi bound under global existence of the toric Calabi flow \cite[Theorems 4 and 16]{Sze08}. 
For a general polarized manifold, long-time solution to the Calabi flow is standing as a difficult probelem and will give a substantial help in proving Theorem \ref{thm:main}, however, we will not take this strategy in this article.

An important step was taken by 
Xia, who proved the corresponding equality on the metric 
completion $(\cE^2,d_2)$ of the space of smooth
K\"ahler potentials $\cH$ \cite[Theorem 1.1(1)]{Xia21}. 
More precisely, he replaced $\Cal(\omega)$ with the maximal rate of decrease  
\begin{equation*}
 \abs{\partial M}(u):=\sup_{v\ne u}\frac{(M(u)-M(v))_+}{d_2(u,v)}
\end{equation*}
of the K-energy at each point $u \in \cE^2$, 
and at the same time replaced $M^{\NA}(\cX, \cL)$ with slope at infinity of 
extended K-energy along finite-energy geodesic rays  (namely, with radial K-energy). 
The remaining two points in proving Theorem \ref{thm:main} are smooth metric approximation: $\inf_{\cH} \Cal(\omega) =\inf_{\cE^2} \abs{\partial M}$  
and approximation of a geodesic ray by test configurations. 

The latter algebraic approximation almost follows from combining recent developments of this area. 
Darvas--Lu prove density of rays with bounded potentials, with simultaneous recovery of radial K-energy \cite[Theorem 1.5(i)]{DL20}. Li proves that finite radial K-energy forces a ray to be maximal \cite[Theorem 1.2]{Li22}. 
He also proves the inequality 
$M^{\NA}(\ell^{\NA})\leq M^{\rad}(\ell)$ 
for an arbitrary maximal ray $\ell$ \cite[Theorem 1.7]{Li22}.  
Now Trusiani's special Fujita approximation \cite[Theorem B and Corollary A]{Tru26} combined with the approximation by model envelopes of Boucksom--Jonsson \cite[Lemmas 7.15--7.17]{BJ25} gives a desired 
sequence of test configurations. 
We need to add one more point to the discussion: convergence in the non-Archimedean $L^1$-distance and
convergence of entropy do not by themselves imply 
convergence of the $L^2$ norms in the denominator of 
\eqref{eq:main}. We obtain the latter convergence serving higher-moment estimate for the solutions
of \eqref{eq:retracted-measures}. 
This completes the proof of 
$\inf_{\cE^2}\abs{\partial M} =\sup \frac{-M^{\NA}(\cX, \cL)}{\norm{(\cX, \cL)}_2}$. 
In addition, While preparing this draft \cite[Theorem 1.5]{DZ26} appeared. It can further 
simplifies our argument. 

The primary focus of this paper is on  
$\inf_{\cH} \Cal(\omega) =\inf_{\cE^2} \abs{\partial M}$ 
and we do not require a polarization at this stage. 
 For $\tau>0$ and $u_0\in\cH$, we consider the  minimizer of  
 \begin{equation}\label{eq:proximal}
 M(u)+\frac{d_2(u,u_0)^2}{2\tau}
 \end{equation}
 on $\cE^2$; see Lemma \ref{lem:resolvent}, 
 and call it {\em K-resolvent}. 
This is known as proximal operator in convex optimization. 
To approximate the resolvent by smooth potentials, we solve 
\eqref{eq:cylinder}--\eqref{eq:endpoint}, which couples a continuous path of prescribed scalar curvature equation  
to a version of Chen's $\e$-subgeodesic; a non-degenerate $(n+1)$-variables complex Monge--Amp\`ere equation on a cylindrical domain \cite[Section 5]{Chen00}. 
The path endpoint at $s=0$ gives the approximate resolvent equation as in \cite[Lemma 4.9]{Streets14}, while $s=1$ easily has an explicit (actually unique) solution. 
Coupling with $\e$-geodesic controls 
 the terminal velocity appearing in the first variation of the proximal objective function \eqref{eq:proximal} and provides 
a uniform in $s$ estimate for the solution 
which enables us to apply Leray-Schauder degree argument. 
The $C^0$-estimate in Lemma \ref{lem:CC} adapts the arguments in Chen--Cheng \cite[Section 5]{CC21I}, where they used a different coupling to solve the cscK equation under the coercivity condition of K-energy. 
In our situation only an upper bound for the function
$f$ in \eqref{eq:CC} is available. 
Lemma \ref{lem:CC} first bounds the endpoint potential $\varphi$. 
Subsequent terminal velocity estimate 
which is equivalent to the boundary mixed second derivative estimate for 
the solution of $(n+1)$-variable Monge--Amp\`ere equation then
gives a lower bound for $f$, so that
\cite[Theorem 5.3]{CC21I} applies. 
After obtaining a solution at the endpoint $s=0$, letting $\e \to 0, \tau \to \infty$ 
completes the proof.

\medskip
\noindent{\bfseries Acknowledgments and AI disclosure.}
The author is supported by JSPS Grant-in-Aid for Scientific Research $(\mathrm{C})$, No.\ 26K06811. 
ChatGPT-6 Astra was used extensively to explore proof candidates and assist with checking calculations, references, and exposition. The author started with a few ideas including avoidance of Calabi flow, use of proximal equation and Trusiani's approximation, but most of other arguments, represented by coupling with $\e$-geodesic, Leray-Schauder degree argument, and use of Boucksom--Jonsson's envelops approximation were suggested by Astra. 
The author directed the research, verified all mathematical arguments and references, supplied the necessary details, and wrote the final exposition. 
The author takes sole responsibility for the contents of this paper. 

%For these reasons, of course, 
%I do not wish to claim originality for much of the mathematical content presented here. 
%Still bewildered by the circumstances in which we now find ourselves, I would like to share the mathematical result we have obtained and discuss it with others. 
%At a time when many long-held values are being called into question, I have been discussing with friends whom I deeply respect and received their thoughtful and candid advice. I wish to emphasize here that what I value most is my connection with people like them. 

\section{Preliminaries on K\"ahler potentials and test configurations}\label{sec:normalization}
\subsection{Energy functionals}
We recall definitions of several functionals on 
the space of K\"ahler potentials 
\eqref{eq:cH}. 
The Monge--Amp\`ere energy 
\begin{equation}\label{eq:energy}
 E(\varphi):=\frac1{(n+1)V}\sum_{j=0}^n
 \int_X\varphi \omega_\varphi^j\wedge\omega_0^{n-j} 
\end{equation}
plays a basic role in the study of general Monge--Amp\`ere equations 
because it satisfies the differential formula 
\begin{equation*}
(dE)_\varphi(f)=\frac1V\int_Xf \omega_\varphi^n. 
\end{equation*}
Twisting by a smooth closed real $(1,1)$-form $\theta$ is also needed. We set
\begin{equation*}
 E_\theta(\varphi)=\frac1V\sum_{j=0}^{n-1}
 \int_X\varphi \theta\wedge\omega_\varphi^j\wedge\omega_0^{n-1-j}.
\end{equation*}
Its differential is $V^{-1}\int_Xf\tr_{\omega_\varphi}\theta \omega_\varphi^n$. With $\bar\theta=n[\theta]\cdot[\omega_0]^{n-1}/V$, put $J_\theta=E_\theta-\bar\theta E$. In particular $J_\theta$ is invariant under the scaling $\varphi \mapsto \varphi+ c$. 

K-energy is characterized by the differential formula \eqref{eq:mabuchi} 
up to an addition of a constant. 
In practice, we use Chen's expression \cite{Chen00Lower} using relative entropy 
\begin{equation*}
 H(\varphi)=\frac1V\int_X\log\frac{\omega_\varphi^n}{\omega_0^n}\omega_\varphi^n.
\end{equation*}
Under these conventions, K-energy is 
written as
\begin{equation}\label{eq:chen-tian}
 M(\varphi)=H(\varphi)+\bar S E(\varphi)-E_{\Ric(\omega_0)}(\varphi). 
\end{equation}
The fundamental feature of K-energy is  convexity along the natural metric structure of $\cH$. 
Following the work of Mabuchi \cite{Mab87} and Darvas, for each $p \geq 1$ we define the $L^p$ length structure by the tangent norm
\begin{equation*}
 \norm{f}_{p,\varphi}=\left[\frac1V\int_X\abs{f}^p\omega_\varphi^n\right]^{\frac{1}{p}}.
\end{equation*}
The associated distance is denoted $d_p$. 
Historically $d_2$ was first studied. 
The metric structure $d_1$ is closely related to the Monge--Amp\`ere energy and plays an important role in variational approach to cscK metrics.  
As it was shown in \cite{LV13}, 
geodesics with smooth endpoints are not necessarily contained in $\cH$. 
To deal with this and to consider the metric completion of $\cH$, we need the Monge--Amp\`ere operator for a possibly singular $\omega_0$-plurisubharmonic potential $u$. 
The non-pluripolar product 
\begin{equation*}
 \MA(u):=V^{-1}\langle(\omega_0+dd^cu)^n\rangle.
\end{equation*}
of Boucksom, Eyssidieux, Guedj, and Zeriahi \cite[Section 1.2]{BEGZ10} is perfectly suited for the purpose. 
We also refer to the textbook \cite{GZ17} 
for basic materials of variational approach to complex Monge--Amp\`ere equations and 
proofs of the facts below. 
For smooth $u$, $\MA(u)$ is the ordinary wedge product $V^{-1}\omega_u^n$. 
In general $\int_X \MA(u) \leq 1 $ and 
$u$ is called full-mass if the equality holds. 
By definition, the finite-energy class $\cE^p$ consists of the full-mass potentials for which $\int_X\abs{u}^p\MA(u)<\infty$. 
If one extends definition of Monge--Amp\`ere energy to arbitrary 
 $\omega_0$-psh functions by 
\begin{equation}\label{eq:extend E}
      E(u) := \inf \bigg\{ E(v): v\in \cH, \ \ u\leq v  \bigg\} \in \R\cup \{-\infty\}, 
\end{equation}
it follows $u \in \cE^1$ if and only if $E(u)>-\infty$. 
The metric structures are 
also controlled by the  Monge-Amp\`ere measures. 

\begin{prop}[Darvas {\cite[Theorem 3]{Dar15}}]\label{prop:dp}
For $p\geq1$, $(\cE^p,d_p)$ is the metric completion of $\cH$. There exists a constant $C_p=C(n,p)>0$ such that, for $u,v\in\cE^p$,
\begin{equation}\label{eq:dp-comparison}
 C_p^{-1}d_p(u,v)^p\leq
 \int_X\abs{u-v}^p(\MA(u)+\MA(v))\leq C_p d_p(u,v)^p.
\end{equation}
Monotone convergence within $\cE^p$ implies convergence in $d_p$.
\end{prop}

The following CAT(0)-property is unique to 
$p=2$. 

\begin{prop}[Darvas {\cite[Theorem 1]{DarCompletion}}]\label{prop:cat0}
The space $(\cE^2,d_2)$ is complete and CAT(0). If $u_t$, $0\leq t\leq1$, is the constant-speed geodesic joining $u_0,u_1\in\cE^2$, then, for any $v\in\cE^2$,
\begin{equation*}
 d_2(v,u_t)^2\leq(1-t)d_2(v,u_0)^2+t d_2(v,u_1)^2
                  -t(1-t)d_2(u_0,u_1)^2.
\end{equation*}
\end{prop}

By observation of Semmes \cite{Sem92}, 
geodesic equation on $\cH$ itself is 
translated into the degenerate Monge--Amp\`ere equation on the 
product $(n+1)$-variable space 
\begin{equation*}
X \times \{ \tau \in \C: e^{-r/2}<\abs{\tau} <1\}. 
\end{equation*}
More precisely, denote a pull-backed semipositive form  $\Omega_0:= p_1^*\omega_0$. 
For a smooth path $u_t \in \cH$ $(0\leq t\leq r)$, the function $U(x, t):= u_{-\log{\abs{\tau}^2}}(x)$ 
on the above product space satisfies 
\begin{equation}\label{eq:geodesic}
(\Omega_0 +dd^c U)^{n+1} \equiv 0 
\end{equation}
if and only if $u_t$ is a geodesic. 
In terms of a smooth path $u_t$, geodesic curvature 
is $\ddot{u}-\abs{\partial \dot{u}}^2_{\omega_u}$. 
Among smooth subsolutions, a smooth solution of the geodesic equation is characterized by the affineness of $E(u_t)$. 
As we have already noted, even for given smooth endpoints $u_0, u_r$, 
the solution of \eqref{eq:geodesic} is no longer smooth. 
The best possible $C^{1,1}$-regularity was obtained 
by \cite{CTW17}. 
More generally, given $u_0, u_1 \in \cE^p$, 
one can consider a connecting geodesic purely 
in terms of metric space.  
The CAT(0)-property assures the uniqueness in $\cE^2$, however, it is not necessarily unique 
in $\cE^1$. 
A special type of $\cE^1$-geodesic arising as a weak solution of \eqref{eq:geodesic} is called 
{\em weak geodesic} or psh geodesic in the 
literature. 
More generally, subsolution to \eqref{eq:geodesic} 
is called {\em subgeodesic}. 
Rigorously, weak geodesic is defined as the envelope of subsolutions \cite[section 5]{DarCompletion}, 
as an analogue of classical Perron-Bremermann method. 
Among all subgeodesics, weak geodesity is still equivalent to  
affineness of $E(u_t)$, if one takes the definition  \eqref{eq:extend E}.

The entropy expression \eqref{eq:chen-tian} extends $M$ to $\cE^1$. Its restriction to $\cE^2$ is lower semicontinuous and geodesically convex. For a positive smooth closed twist $\chi$, the same holds for $M+sJ_\chi$, $s\geq0$; see \cite[Theorem 1.2]{BDL17} which extends the weak-geodesic convexity of Berman--Berndtsson \cite[Theorem 1.1]{BB17}. 
 More generally, for a proper lower semicontinuous convex functional $F$ on a geodesic metric space $(M,d)$, its maximal rate of decrease at $u\in\{F<\infty\}$ is
\begin{equation*}
 \abs{\partial F}(u):=\sup_{v\ne u}\frac{(F(u)-F(v))_+}{d(u,v)}. 
\end{equation*}
We put $\abs{\partial F}(u):=+\infty$ when $F(u)=+\infty$. At a smooth potential the K-energy slope satisfies
\begin{equation}\label{eq:smooth-slope}
 \abs{\partial M}(u)=\Cal(\omega_u)
\end{equation}
by \cite[Proposition 3.6]{Xia21}. 
The left-hand side therefore 
gives a natural extension of Calabi functional to $\cE^2$. 
Note that this identification uses the smooth Calabi flow but only on a short time. 
We refer \cite{AGS08} for general theory 
of convex optimization in a metric space. 

We end this subsection raising two propositions used later. 

\begin{prop}[Berman--Darvas--Lu {\cite[Lemma 3.3]{BDL19}}]\label{prop:mabuchi-support}
Let $u_0\in\cH$, $u_1\in\cE^1$, and let $u_t$, $0\leq t\leq1$, be the finite-energy weak geodesic segment with these endpoints. For a smooth positive closed $(1,1)$-form $\chi$ and $s\geq0$, 
\begin{equation*}
 \begin{split}
 (M+sJ_\chi)(u_1)-(M+sJ_\chi)(u_0)
 \geq-\frac1V\int_X[S_{u_0}-\bar S-s(\tr_{\omega_{u_0}}\chi-\bar\chi)]\dot u_0 \omega_{u_0}^n.
 \end{split}
\end{equation*}
Here $S_u:=S(\omega_u)$ and $\dot u_0$ is the right derivative at the initial endpoint of this segment. 
\end{prop}

\begin{proof}
For $s>0$, apply the supporting-slope inequality of
\cite[Lemma 3.3]{BDL19} to the positive form $s\chi$.
After division by the fixed volume, their twisted K-energy is $M+sJ_\chi$, up to an additive constant. 
Their Aubin--Mabuchi energy is our $E$. 
Put $m(t):=(M+sJ_\chi)(u_t)$. Convexity gives
\begin{equation*}
 m(1)-m(0)\geq\frac{m(t)-m(0)}t\qquad(0<t\leq1).
\end{equation*}
The cited lemma identifies a lower bound for the lower limit of the
quotient as $t \to 0$:
\begin{equation*}
 \liminf_{t\to 0}\frac{m(t)-m(0)}t
 \geq-\frac1V\int_X
 [S_{u_0}-\bar S-s(\tr_{\omega_{u_0}}\chi-\bar\chi)]
 \dot u_0 \omega_{u_0}^n.
\end{equation*}
Combining these inequalities proves the assertion. s
If $M(u_1)=+\infty$, the assertion is immediate. Otherwise, letting $s\to 0$ gives the untwisted case. 
\end{proof}

\begin{prop}[Berman--Darvas--Lu {\cite[Theorem 1.3]{BDL17}}]\label{prop:entropy-recovery}
For any $u\in\cE^2$ with $M(u)<\infty$, there exists a sequence $u_j\in\cH$ such that 
\begin{equation*}
 d_2(u_j,u)\longrightarrow0,\qquad H(u_j)\longrightarrow H(u), 
 \qquad M(u_j)\longrightarrow M(u).
\end{equation*}
If $E(u)=0$, the approximants can be chosen with $E(u_j)=0$.
\end{prop}
\begin{proof}
Apply the cited theorem with $p=2$ and reference volume $\omega_0^n$. The nonentropy terms in \eqref{eq:chen-tian} are $d_1$-continuous. Subtracting $E(u_j)$ preserves entropy and K-energy, and this constant converges to $E(u)$.
\end{proof}

\subsection{Test configurations and the nonreduced-fiber correction}\label{sec:test configurations}
We basically follow the terminology in \cite{BHJ17}, \cite{BHJ19}. 
A test configuration of exponent $r$ consists of a normal variety $\cX$, a flat projective map $\pi:\cX\to\C$, and a relatively ample line bundle $\cA$ endowed with a $\C^*$-action lifting multiplication on the base. Over $\C^*$ it is equivariantly isomorphic to the product of $(X,L^r)$ and $\C^*$. We write $\Q$-line bundle $\cL=r^{-1}\cA$. Gluing to the trivial family near infinity, as in \cite[Section 2.2]{BHJ17}, gives the compactification $(\overline{\cX},\overline{\cL})\to\mP^1$.

Let us for a while fix a smooth $\mathbb{S}^1$-invariant Hermitian form 
$\Omega$ in $c_1(\overline{\cL})$. 
equivariant trivialization identifies it with a family of metric 
$\omega_0+dd^c\psi_t \in c_1(L)$, with $t=-\log{\abs{\tau}^2}$. 
A weak geodesic ray $\ell$ in $\cE^1$ with
$\ell_0\in\cH$ 
is {\em associated with}, or compatible with $(\cX, \cL)$ 
if it satisfies 
$\abs{\ell_t -\psi_t} \leq C$ uniformly in $t$.  
Equivalently, the potential on the punctured product space, induced by the ray, extends across the central fiber as a locally bounded plurisubharmonic weight of $\cL$. 
The associated ray exists and is unique for the prescribed
initial potential. See \cite[Corollary 5.3]{BJ25}.
For this reason test configurations can be regarded as 
algebraic counterparts of weak geodesic rays. 
Furthermore, along the associated weak geodesic ray, 
slope of K-energy at infinity is precisely equivalent to 
non-Archimedean K-energy $M^\NA(\cX, \cL)$ whose precise definition will be given below \cite[Theorem C]{BJ25}. 
This fact already has its roots in the results of Ding--Tian and supports legitimacy of 
Yau--Tian--Donaldson conjecture. 

If $\lambda_j^{(k)}$ are the weights on $H^0(\cX_0,\cA_0^k)$, set $N_k:=h^0(X,L^{rk})$. The Duistermaat--Heckman measure is
\begin{equation*}
 \mu_{(\cX,\cL)}:=\lim_{k\to\infty}\frac1{N_k}
             \sum_{j=1}^{N_k}\delta_{\frac{\lambda_j^{(k)}}{rk}}.
\end{equation*}
The limit is weak convergence of probability measures; see \cite[Definition 3.10]{BHJ17}. 
The $L^2$-norm of a test configuration 
is by definition the square root of the variance of $\mu_{(\cX, \cL)}$. 
Twisting the $\C^*$-action by a character translates the weights and leaves \eqref{eq:tc-norm} unchanged. 

\begin{prop}[{\cite[Theorems 1.1 and 1.2]{His16}}]\label{prop:dh-measure}
Let $\ell_t$ be the weak geodesic ray associated to $(\cX,\cL)$. The push-forward by velocity measure is independent of $t$ and satisfies
\begin{equation}\label{eq:dh-measure}
 (\dot\ell_t)_*\MA(\ell_t)=\mu_{(\cX,\cL)}.
\end{equation}
Consequently the centered ray $\ell_t-E(\ell_t)+E(\ell_0)$ has $d_2$ speed $\norm{(\cX,\cL)}_2$.
\end{prop}

Following \cite[Definition 7.13]{BHJ17}, 
we define the Donaldson--Futaki invariant and the non-Archimedean K-energy 
by the intersection numbers 
\begin{align}
 &\DF(\cX,\cL)
 =\frac{K_{\overline{\cX}/\mP^1}\overline{\cL}^{n}}V
 +\frac{\bar S \overline{\cL}^{n+1}}{(n+1)V},\nonumber\\
 &M^{\NA}(\cX,\cL)
 =\DF(\cX,\cL)+\frac{(\cX_{0,\red}-\cX_0)\cL^n}V.
 \label{eq:correction}
\end{align}
In particular $\DF \geq M^{\NA}$. Equality holds when the central fiber is reduced. 
Under normalized base change of degree $d$, $M^{\NA}$ and the norm both multiply by $d$ \cite[Lemma 6.9 and Proposition 7.14]{BHJ17}, 
whereas $\DF$ need not scale homogeneously. 
In our situation, however, 
either invariant can be used for the formulation. 

\begin{lem}\label{lem:base-change}
The supremum in \eqref{eq:main} is unchanged if $M^{\NA}$ is replaced by $\DF$, or if it is restricted to reduced central fibers.
\end{lem}
\begin{proof}
Take the normalized base change of degree $d$ of a given test configuration. 
Normalization preserves the Duistermaat--Heckman measure \cite[Theorem 3.14]{BHJ17}. 
Relative ampleness is also  preserved. 
For sufficiently divisible $d$, the new central fiber is reduced and its Donaldson--Futaki invariant is $dM^{\NA}(\cX,\cL)$ \cite[Proposition 7.16]{BHJ17}.  Thus each ratio for $M^{\NA}$ is realized by a ratio for $\DF$ with reduced central fiber. Conversely, \eqref{eq:correction} bounds every $-\DF$ ratio above by the corresponding $-M^{\NA}$ ratio. Restriction to reduced fibers cannot increase the supremum. 
\end{proof}

\subsection{Geodesic rays and Xia's formulation}
Let $\cR^2$ consist of constant-speed weak geodesic rays $\ell:[0,\infty)\to\cE^2$ based at a fixed smooth potential. Define radial K-energy as 
\begin{equation*}
 M^{\rad}(\ell):=\lim_{t\to\infty}\frac{M(\ell_t)}t.
\end{equation*}
Convexity ensures existence of the limit in $\R \cup \{ \infty\}$. 
Now we can state Xia's result. 

\begin{prop}[Xia {\cite[Theorem 1.1(1)]{Xia21}}]\label{prop:xia}
\begin{equation}\label{eq:xia}
 \inf_{u\in\cE^2}\abs{\partial M}(u)=\max_{\ell\in\cR^2\setminus\{0\}}
 \frac{-M^{\rad}(\ell)}{d_2(\ell_0, \ell_1)}.
\end{equation}
\end{prop}
The inclusion of constant-translation rays makes this maximum nonnegative. 
When the left-hand side is positive, this maximum is attained by a geodesic
ray with finite radial K-energy
\cite[Theorem 1.1(1)]{Xia21}. 
For a weak geodesic ray associated with a test configuration, $M^{\rad}\leq M^{\NA}$ \cite[Proposition 5.1]{Xia21}. Centering its velocity, applying \eqref{eq:xia}, and then taking the supremum proves $\inf_{\cE^2}\abs{\partial M} \geq \frac{-M^{\NA}(\cX, \cL)}{\norm{(\cX, \cL)}_2}$. On the other hand, formula~\eqref{eq:smooth-slope} and $\cH\subset\cE^2$ show $\inf_{\cH} \Cal \geq \inf_{\cE^2}\abs{\partial M}$.

\section{Proximal operators and the infimum of the Calabi functional}\label{sec:resolvent}
We first recall definition of the proximal operator $p_{\tau}$ for a convex
functional on a complete CAT(0) space. 
See \cite{Bac14} for a comprehensive treatment. 
The following lemma shows that
$d(u,p_\tau)/\tau$ converges to the infimum of the
maximal rate of decrease as $\tau\to\infty$. 
For a Hilbert space it is intuitively clear because the minimizer enjoys $\nabla F + (p_\tau -u)/\tau=0$. 
Convexity qualifies the general case. 

\begin{lem}[See also {\cite[section 3.1]{AGS08}}]\label{lem:resolvent}
Let $(M,d)$ be complete and CAT(0), and let $F$ be proper, lower semicontinuous, and geodesically convex. Suppose some point has finite descending slope. For fixed $u\in M$, the unique minimizer
\begin{equation*}
 p_\tau=\argmin_{w\in M}\left(F(w)+\frac{d(u,w)^2}{2\tau}\right)
\end{equation*}
exists for every $\tau>0$. If $m_F:=\inf_M\abs{\partial F}$, then
\begin{equation}\label{eq:resolvent}
m_F\leq\abs{\partial F}(p_\tau)
 \leq\frac{d(u,p_\tau)}\tau
 \leq\abs{\partial F}(v)+\frac{2d(u,v)}\tau
\end{equation}
for every finite-slope point $v$. In particular $d(u,p_\tau)/\tau\to m_F$ as $\tau\to\infty$.
\end{lem}
\begin{proof}
The functional to be minimized: 
\begin{equation*}
F_\tau(w):= F(w)+\frac{d(u, w)^2}{2\tau} 
\end{equation*}
is coercive. Namely, There exists constants $\e, C>0$ such that $F_\tau \geq \e d(u,\cdot)-C$.  
It follows from the slope inequality $F(w)\geq F(v)-\abs{\partial F}(v)d(v,w)$. 
The CAT(0) inequality in Proposition \ref{prop:cat0} implies that 
$F_\tau$ is $1/\tau$-strongly convex. 
In particular, the midpoint $w_{\frac{1}{2}}$ of $w_0,w_1$, the inequality reads $d(u,w_{\frac{1}{2}})^2\leq\tfrac12d(u,w_0)^2+\tfrac12d(u,w_1)^2-\tfrac14d(w_0,w_1)^2$ so that 
\begin{equation*}
F_\tau(w_\frac{1}{2})
\leq \frac{F_\tau(w_0)+F_\tau(w_1)}{2} 
-\frac{d(w_0, w_1)^2}{8\tau}. 
\end{equation*}
It implies that every minimizing sequence is Cauchy. Completeness of $M$ and lower semicontinuity of $F$ give existence of the minimizer. Strong convexity gives uniqueness.

Applying minimality of $p_\tau$ and strong convexity on the segment from $p_\tau$ to $v$ gives
\begin{equation*}
 -\abs{\partial F}(v)d(v, p_\tau)\leq F(p_\tau)-F(v)
 \leq\frac{d(u, v)^2-d(u,p_\tau)^2-d(p_\tau, v)^2}{2\tau}.
\end{equation*}
Since $d(u,p_\tau)\geq\abs{d(p_\tau, v)-d(u, v)}$, this implies $d(p_\tau, v)\leq d(u, v)+\tau \abs{\partial F}(v)$. Triangle inequality then gives the last inequality in \eqref{eq:resolvent}. 

For the middle inequality in \eqref{eq:resolvent}, let $w_t$ be the segment from $p_\tau$ to $w$. Minimality and convexity give
\begin{equation*}
 t(F(p_\tau)-F(w))
 \leq\frac{d(u,w_t)^2-d(u,p_\tau)^2}{2\tau}
 \leq\frac{td(u,p_\tau)d(p_\tau,w)}\tau
       +\frac{t^2d(p_\tau,w)^2}{2\tau}.
\end{equation*}
The last inequality is the triangle inequality squared. If we divide by $t$, let $t$ decrease to zero, and take the supremum over $w$, it then completes the proof. 
\end{proof}

Let us apply this to $M=\cE^2$, K-energy, and a fixed $u_0\in\cH$ with $E(u_0)=0$. Smooth points have finite slope by \eqref{eq:smooth-slope}. To prove $\inf_{\cH}\Cal\leq\inf_{\cE^2}\abs{\partial M}$, it is enough to construct, for each $\tau>0$, smooth $\varphi_{\tau,\e}$ such that
\begin{equation*}
 \limsup_{\e \to 0}\Cal(\omega_{\varphi_{\tau,\e}})
 \leq\frac{d_2(u_0,p_\tau)}\tau.
\end{equation*}
The rest of the analytic proof is devoted to constructing these smooth metrics. 

\section{A coupled elliptic system and energy comparison}\label{sec:coupled}
The boundary $t=1$ equation below is motivated by the first variation of the proximal objective function. Streets already writes the functional $M+d_2(u_0,\cdot)^2/(2\tau)$ in \cite[Definition 4.4]{Streets14} and computes its first variation at a smooth endpoint in \cite[Lemma 4.9]{Streets14}. The $\e$-geodesic is a positive-density regularization of the geodesic equation on $\cH$, as in \cite[Section 5, Definition 4]{Chen00}. 
We choose the numerator in the right-hand side of \eqref{eq:cylinder} which gives the geodesic curvature, to be the volume form of the endpoint $s=1$ of the continuity path. This couples the Dirichlet problem to a prescribed scalar curvature equation. 

Fix $u_0\in\cH$ with $E(u_0)=0$ and a parameter $\tau>0$. If $\lambda_{\min}(x)$ is the pointwise least eigenvalue of $\Ric(\omega_{u_0})$ relative to $\omega_{u_0}$, choose
\begin{equation}\label{eq:twist-choice}
 \Lambda:=1+\max\{0,-\min_X\lambda_{\min}\},\qquad
 \chi:=\Ric(\omega_{u_0})+\Lambda\omega_{u_0}>0.
\end{equation}
Set $M_s:=M+sJ_\chi$ for $0\leq s\leq1$, and fix $0<\e<\min\{1,\tau \Lambda\}$. The parameter $s$ will be used for the continuity path. The parameter $t\in[0,1]$ is for a path $\gamma:[0,1]\to\cH$ from $u_0$ to an unknown potential $\varphi$. A dot like $\dot{\gamma}$ always denotes differentiation in $t$. We consider
\begin{align}
 &\ddot\gamma-\abs{\partial\dot\gamma}_{\omega_{\gamma(t)}}^2
   =\e\frac{\omega_\varphi^n}{\omega_{\gamma(t)}^n},
   \qquad \gamma(0)=u_0,\quad\gamma(1)=\varphi,\label{eq:cylinder}\\
 &S_\varphi-\bar S-s(\tr_{\omega_\varphi}\chi-\bar\chi)
   =\frac{\dot\gamma(1)}\tau.\label{eq:endpoint}
\end{align}
For $\e=0$, 
\eqref{eq:cylinder} is precisely 
the geodesic equation \eqref{eq:geodesic}. 
For $\e=0$ and $s=0$, let $\gamma$ be such a weak geodesic from $u_0$ to $\varphi$. The first-variation formula of Streets gives, for a smooth variation $f$ of the endpoint,
\begin{equation*}
 d\left(M+\frac{d_2(u_0,\cdot)^2}{2\tau}\right)_\varphi(f)
 =-\frac1V\int_X f\left(S_\varphi-\bar S-
                    \frac{\dot\gamma(1)}\tau\right)\omega_\varphi^n.
\end{equation*}
Thus \eqref{eq:endpoint} at $\e=s=0$ 
gives a critical point for the functional $F_\tau $ in Lemma \ref{lem:resolvent}. 

For fixed smooth $u_0,\varphi$ and $\e>0$, equation
\eqref{eq:cylinder} has a unique smooth solution.
The corresponding nondegenerate Dirichlet problems on a 
cylindrical domain also originate from 
\cite{Chen00}. Subsolution method on K\"ahler manifolds is discussed in \cite{Guan09,Blocki12}.
For the precise existence theorem used here, we apply
\cite[Theorem 1.1]{GL10}.
Its hypotheses and the local dependence on the endpoint are
verified in Proposition~\ref{prop:dirichlet}.
Until then we work with an a priori smooth solution of equations \eqref{eq:cylinder}-\eqref{eq:endpoint}.

\subsection{Two lemmas concerning velocities and geodesic curvatures}
This subsection is preparation to section \ref{sec:energy comparison}. 
The following local barrier calculation gives the second-order endpoint expansion directly from the weak-geodesic equation.

\begin{lem}\label{lem:quadratic-approximation}
Let $\gamma\colon [0,1]\to\cH$ be a smooth path with uniformly positive spatial metrics. Write $c(\gamma):=\ddot\gamma-\abs{\partial\dot\gamma}^2_{\omega_\gamma}$ for its geodesic curvature. For a real increment $h$ with $t,t+h\in[0,1]$, let $\xi(t,h)$ denote the initial velocity of the weak geodesic from $\gamma(t)$ to $\gamma(t+h)$. Then it holds 
\begin{equation}\label{eq:quadratic-approximation}
 \norm{\xi(t,h)-h\dot\gamma(t)-\tfrac12h^2c(\gamma)(t)}_\infty
 \leq A\abs{h}^3
\end{equation}
for sufficiently small $\abs{h}$, with $A$ independent of $t$.
\end{lem}
\begin{proof}
Set $\delta:=\gamma(t+h)-\gamma(t)$ and
$B_\delta:=\abs{\partial\delta}_{\omega_{\gamma(t)}}^2$.
For $s\in[0,1]$, we denote 
\begin{equation*}
 q(s):=\gamma(t)+s\delta+\tfrac12(s^2-s)B_\delta.
\end{equation*}
Its spatial metric is positive for small $\abs{h}$. Indeed,
$\norm{\delta}_{C^4}=O(\abs{h})$ and $\norm{B_\delta}_{C^2}=O(h^2)$
uniformly along the fixed smooth path. Differentiation with respect to $s$ gives
\begin{equation*}
 c(q):=\partial_s^2q-\abs{\partial\partial_sq}_{\omega_q}^2
 =B_\delta-\abs{\partial\delta+(s-1/2)\partial B_\delta}_{\omega_q}^2
 =O(\abs{h}^3).
\end{equation*}
Here changing the inverse spatial metric contributes $O(\abs{h}^3)$,
as does the cross term containing $\partial B_\delta$; its square is $O(h^4)$.
Choose $A_0$ so that $\abs{c(q)}<2A_0\abs{h}^3$, and 
define the two barriers as 
\begin{equation*}
 q_\pm(s):=q(s)\pm A_0\abs{h}^3s(1-s). 
\end{equation*}
These additions are spatially constant. Consequently
\begin{equation*}
 c(q_-)=c(q)+2A_0\abs{h}^3>0,\qquad
 c(q_+)=c(q)-2A_0\abs{h}^3<0.
\end{equation*}
Let $\alpha(s)$ be the weak geodesic with the same endpoints.
The subgeodesic $q_-$ is a competitor in its envelope definition, so
$q_-\leq\alpha$. To prove the opposite barrier inequality, suppose
$m:=\max(\alpha-q_+)>0$. Since the endpoint values agree, the maximum
is attained at an interior point. At that point $q_++m$ touches
$\alpha$ from above. After adding a local potential for $\omega_0$,
this is a smooth upper test for a plurisubharmonic function; its $(n+1)$-variable 
complex Hessian must therefore be positive semidefinite. Its spatial
block is positive definite, but its Schur complement has the sign of
$c(q_+)<0$. This is 
a contradiction. Hence
\begin{equation*}
 q_-\leq\alpha\leq q_+,\qquad
 \abs{\alpha(s)-q(s)}\leq A_0\abs{h}^3s(1-s).
\end{equation*}
Subtracting the common initial value and dividing by $s$ gives
$\xi(t,h)=\delta-B_\delta/2+O(\abs{h}^3)$ at $s=0$.
Taylor's theorem yields
\begin{equation*}
 \delta=h\dot\gamma(t)+\tfrac12h^2\ddot\gamma(t)+O(\abs{h}^3),
 \qquad B_\delta=h^2\abs{\partial\dot\gamma(t)}_{\omega_{\gamma(t)}}^2
 +O(\abs{h}^3),
\end{equation*}
with remainders in $C^1$ and $C^0$, respectively. Substitution proves
\eqref{eq:quadratic-approximation}.
\end{proof}

For $u,v\in\cH$, let $u_t$ be the weak geodesic segment from $u$ to $v$, and define $\xi_u(v):=\dot u_0$. Its endpoint speed satisfies
\begin{equation}\label{eq:endpoint-speed}
 \norm{\xi_u(v)}_{2,u}=d_2(u,v).
\end{equation}
This is the geodesic speed formula recalled in \cite[Section 2.1, equation (2.1)]{Xia21}.

\begin{prop}\label{prop:distance-support}
For smooth potentials $u,v,u',v' \in \cH$, one has
\begin{equation}\label{eq:distance-support}
 \frac{d_2(u',v')^2-d_2(u,v)^2}{2}
 \geq-\langle\xi_u(v),\xi_u(u')\rangle_u
       -\langle\xi_v(u),\xi_v(v')\rangle_v,
\end{equation}
where the bracket is the $L^2(V^{-1}\omega_u^n)$ pairing at $u$. 
\end{prop}
\begin{proof} 
We first compute the right derivative of square distance along geodesics. 
Weak geodesics with smooth endpoints have $C^{1,1}$ regularity by \cite[Theorem 1.1]{CTW17}. Their endpoint speed is given by \eqref{eq:endpoint-speed}. For two such geodesics $\alpha,\beta$ starting at $u$ with velocities $\xi,\eta$, this gives
\begin{equation*}
 \alpha(t)=u+t\xi+O(t^2),\qquad \beta(t)=u+t\eta+O(t^2)
\end{equation*}
in $C^0$, with continuous velocities. We approximate $\xi,\eta$ uniformly by fixed smooth functions $f,g$. Lemma \ref{lem:quadratic-approximation} and \eqref{eq:endpoint-speed} then give
\begin{equation*}
 \lim_{t\to 0}\frac{d_2(u+tf,u+tg)}t
 =\norm{f-g}_{2,u}.
\end{equation*}
The estimate $d_2(v,w)\leq C\norm{v-w}_\infty$ for bounded potentials follows from \eqref{eq:dp-comparison}. It permits passage first to $t\to 0$, and then to the uniform smooth-approximation limit. Consequently $d_2(\alpha(t),\beta(t))/t\to\norm{\xi-\eta}_{2,u}$. When both velocities are nonzero, writing $\angle_u(\alpha,\beta)$ for the Alexandrov angle, the Euclidean comparison cosine law gives
\begin{equation*}
 \cos\angle_u(\alpha,\beta)=\lim_{t\to 0}
 \frac{d_2(u,\alpha(t))^2+d_2(u,\beta(t))^2-d_2(\alpha(t),\beta(t))^2}
 {2d_2(u,\alpha(t))d_2(u,\beta(t))}
 =\frac{\langle\xi,\eta\rangle_u}{\norm{\xi}_{2,u}\norm{\eta}_{2,u}}.
\end{equation*}
The CAT(0) first-variation formula \cite[II.3, Corollary 3.6]{BrHa99} hence gives the right derivative
\begin{equation*}
 \left.\frac{d}{dt}\right|_{t=0+}\frac12d_2(\alpha(t),v)^2
 =-\langle\xi_u(v),\xi\rangle_u.
\end{equation*}
If $\xi=0$, the geodesic $\alpha$ is constant and both sides vanish. If $v=u$, then $d_2(\alpha(t),u)^2=t^2\norm{\xi}_{2,u}^2$, so both sides again vanish. These cases require no angle at a zero velocity.

Now let $u_t$ and $v_t$, $0\leq t\leq1$, be the geodesics from $u$ to $u'$ and from $v$ to $v'$.
We claim that \eqref{eq:distance-support} follows from the convexity of $F(t):=d_2(u_t,v_t)^2/2$. 
Indeed, the CAT(0) quadrilateral inequality, equivalently the Cauchy--Schwarz inequality for quasilinearization in \cite{BN08}, gives
\begin{equation*}
 \abs{d_2(u_t,v_t)^2-d_2(u_t,v)^2-d_2(u,v_t)^2+d_2(u,v)^2}
 \leq2t^2d_2(u,u')d_2(v,v').
\end{equation*}
 Subtracting $F(0)$ and dividing by $t>0$ yields
\begin{align*}
 \frac{F(t)-F(0)}t
 &=\frac{d_2(u_t,v)^2-d_2(u,v)^2}{2t}
   +\frac{d_2(u,v_t)^2-d_2(u,v)^2}{2t}+O(t).
\end{align*}
The right-derivative formulas at the two endpoints therefore give
\begin{equation*}
 F'_+(0)=-\langle\xi_u(v),\xi_u(u')\rangle_u
          -\langle\xi_v(u),\xi_v(v')\rangle_v.
\end{equation*}
Convexity of $F$ is 
\begin{equation*}
(1-t)F(0)+tF(1) \geq  F(t),\qquad
 F(1)-F(0) \geq F'_+(0).
\end{equation*}
Substitution proves \eqref{eq:distance-support}.

\end{proof}

\subsection{Energy comparison}\label{sec:energy comparison}
In this and the next subsection, 
we fix a solution $(\gamma, \varphi)$ of 
\eqref{eq:cylinder}-\eqref{eq:endpoint}. 
For any smooth path
$\beta \colon [0,1]\to\cH$ with $\beta(0)=u_0, \beta(1)=\psi$, define the $\e$-kinetic energy as 
\begin{equation*}
 K_{\varphi,\e}(\beta):=\int_0^1\left(
 \frac1{2V}\int_X\abs{\dot\beta(t)}^2\omega_{\beta(t)}^n
 +\frac\e V\int_X\beta(t)\omega_\varphi^n\right)dt.
\end{equation*}
See also the first identity in the 
proof of Proposition \ref{prop:energy-bounds}. 
The density $V^{-1}\omega_\varphi^n$ in the second term is fixed. 
It is not $V^{-1}\omega_\psi^n$. 
The following energy comparison 
is natural in view of proximal equation and plays a key role in proving Theorem \ref{thm:main}. 

\begin{lem}\label{lem:energy-comparison}
Let $(\varphi,\gamma)$ solve \eqref{eq:cylinder}--\eqref{eq:endpoint}. For every $\psi\in\cH$ and every smooth path $\beta$ with $\beta(0)=u_0$, $\beta(1)=\psi$, we have 
\begin{equation}\label{eq:energy-comparison}
 \tau M_s(\varphi)+K_{\varphi,\e}(\gamma)
 \leq\tau M_s(\psi)+K_{\varphi,\e}(\beta).
\end{equation}
\end{lem}
\begin{proof}
Put $h=1/N$ and $t_i=ih$ for $i=0,\ldots,N$. For $i=0,\ldots,N$,
let $\xi_i:=\xi_{\gamma(t_i)}(\beta(t_i))$. Define the discrete
kinetic energy of a path $\beta$ by
\begin{equation*}
 K_N(\beta):=\frac1{2h}\sum_{i=0}^{N-1}
 d_2(\beta(t_i),\beta(t_{i+1}))^2.
\end{equation*}
Applying Proposition~\ref{prop:distance-support} to each pair of
successive vertices and collecting them gives
\begin{align*}
 K_N(\beta)-K_N(\gamma)
 &\geq-\frac1h\sum_{i=1}^{N-1}
 \left\langle\xi_{\gamma(t_i)}(\gamma(t_{i-1}))+
 \xi_{\gamma(t_i)}(\gamma(t_{i+1})),\xi_i\right\rangle_{\gamma(t_i)}\\
 &\quad-\frac1h\langle\xi_\varphi(\gamma(1-h)),\xi_N\rangle_\varphi.
\end{align*}
The initial term is zero because $\xi_0=0$. Lemma \ref{lem:quadratic-approximation}
therefore gives
\begin{equation*}
 K_N(\beta)-K_N(\gamma)
 \geq-h\sum_{i=1}^{N-1}\langle c(\gamma)(t_i),\xi_i\rangle_{\gamma(t_i)}
 +\langle\dot\gamma(1),\xi_N\rangle_\varphi-O(h).
\end{equation*}
We claim that the error is uniform in $N$ and $i$: the endpoints of the comparison geodesics range
in two fixed smooth paths with uniformly positive metrics. 
For each $t$, let $\alpha_t(s)$ $(0\leq s\leq1)$ be the weak
geodesic from $\gamma(t)$ to $\beta(t)$, and set
\begin{equation*}
 q_t(s):=(1-s)\gamma(t)+s\beta(t)-As(1-s),
\end{equation*}
Similarly to the proof of Lemma \ref{lem:quadratic-approximation}, 
$q_t$ is a subgeodesic for all sufficiently large $A$.
Since $q_t(0)=\gamma(t)$ and $q_t(1)=\beta(t)$, the subgeodesic
$q_t$ is a competitor in the envelope defining $\alpha_t$, so 
$q_t(s)\leq\alpha_t(s)$.
Together with pointwise convexity in $s$, this gives
\begin{equation*}
 (1-s)\gamma(t)+s\beta(t)-As(1-s)
 \leq\alpha_t(s)
 \leq(1-s)\gamma(t)+s\beta(t).
\end{equation*}
Subtracting $\gamma(t)$, dividing by $s>0$, and letting
$s\to 0$, we obtain
\begin{equation*}
 \beta(t)-\gamma(t)-A
 \leq \partial_s^+\alpha_t(0)
 \leq \beta(t)-\gamma(t).
\end{equation*}
At $t=t_i$, the middle term is $\xi_i$.
Since $\beta$ and $\gamma$ are fixed smooth paths,
this gives a bound $\norm{\xi_i}_\infty\leq B$
independent of $N$ and $i$.
Pairing with the remainder term 
\begin{equation*}
 \xi_{\gamma(t_i)}(\gamma(t_i-h))
 +\xi_{\gamma(t_i)}(\gamma(t_i+h))
 -h^2c(\gamma)(t_i) 
\end{equation*}
yields 
uniformity of $O(h)$, 
as we claimed. 

Let us return to the proof of the lemma. 
Pointwise convexity gives $\beta(t_i)-\gamma(t_i)\geq\xi_i$.
The equation \eqref{eq:cylinder} gives the equality of measures
\begin{equation*}
 c(\gamma)(t_i)\omega_{\gamma(t_i)}^n=\e\omega_\varphi^n.
\end{equation*}
Consequently the source terms cancel the interior terms as follows:
\begin{align*}
 &K_N(\beta)-K_N(\gamma)
 +\frac{\e h}{V}\sum_{i=1}^{N-1}
     \int_X(\beta(t_i)-\gamma(t_i))\omega_\varphi^n\\
 &\quad\geq
 -h\sum_{i=1}^{N-1}\frac1V\int_X
      c(\gamma)(t_i)\xi_i\omega_{\gamma(t_i)}^n
 +\frac{\e h}{V}\sum_{i=1}^{N-1}\int_X\xi_i\omega_\varphi^n
 +\langle\dot\gamma(1),\xi_N\rangle_\varphi-O(h)\\
 &\quad=\langle\dot\gamma(1),\xi_N\rangle_\varphi-O(h).
\end{align*}
Proposition~\ref{prop:mabuchi-support}, applied from $\varphi$ to $\psi$,
and the endpoint equation yield
\begin{align*}
 \tau(M_s(\psi)-M_s(\varphi))
 &\geq-\frac\tau V\int_X
 [S_\varphi-\bar S-s(\tr_{\omega_\varphi}\chi-\bar\chi)]
 \xi_N\omega_\varphi^n\\
 &=-\langle\dot\gamma(1),\xi_N\rangle_\varphi.
\end{align*}
Adding the last two inequalities cancels the terminal term. For each
fixed smooth path $\beta$, Lemma \ref{lem:quadratic-approximation} and
\eqref{eq:endpoint-speed} imply
\begin{equation*}
 d_2(\beta(t_i),\beta(t_i+h))^2
 =h^2\norm{\dot\beta(t_i)}_{2,\beta(t_i)}^2+O(h^3).
\end{equation*}
Hence $K_N(\beta)$ converges to its kinetic energy. Letting $N\to\infty$
proves \eqref{eq:energy-comparison}.
\end{proof}

\subsection{Uniform bounds for the kinetic energy and entropy}
In this subsection we make use of the energy comparison so that obtain energy bounds along the solution to \eqref{eq:cylinder}-\eqref{eq:endpoint}. 
For $\varphi,\psi\in\cH$, Aubin's functionals are defined as 
\begin{align*}
 I_\psi(\varphi)&:=\frac1V\int_X(\varphi-\psi)   (\omega_\psi^n-\omega_\varphi^n),\\
 J_\psi(\varphi)&:=\frac1V\int_X(\varphi-\psi)\omega_\psi^n
                                       -E(\varphi)+E(\psi).
\end{align*}
In particular we need non-negative $I-J$ functional 
\begin{equation}\label{eq:I-J}
 (I-J)_\psi(\varphi):=E(\varphi)-E(\psi)
                     +\frac1V\int_X(\psi-\varphi)\omega_\varphi^n \ \geq0.
\end{equation}
Note that if a path was rescaled to $E(\beta(t))=0$, 
the second term of the $\e$-kinetic energy is essentially $(I-J)_{\beta(t)}(\varphi)$. 
 The conventions agree with \cite[Section 1.4]{BHJ19}. Integration by parts gives the variation
\begin{equation*}
 d_\varphi(I-J)_{u_0}(\varphi)(f)
 =\frac1V\int_X f\Delta_\varphi(u_0-\varphi)\omega_\varphi^n
 =\frac1V\int_X f(\tr_{\omega_\varphi}\omega_{u_0}-n)\omega_\varphi^n.
\end{equation*}
Consequently $(I-J)_{u_0}(\varphi)=J_{\omega_{u_0}}(\varphi)$ when the twisted functional is normalized to vanish at $u_0$. 
Notice that this twisted functional is different from Aubin's $J_{u_0}$. All the expressions in \eqref{eq:I-J} are unchanged by scaling $\varphi \mapsto \varphi+c, \psi \mapsto \psi+c$. 

Before obtaining uniform bound for the kinetic energy, we shall rescale $\gamma(t)$. 
Integrating \eqref{eq:endpoint} and differentiating \eqref{eq:energy} twice along \eqref{eq:cylinder} give
\begin{equation*}
 \frac{d}{dt}E(\gamma(t))\bigg|_{t=1}=0,\qquad
 \frac{d^2}{dt^2}E(\gamma(t))=\e,\qquad E(\gamma(0))=0.
\end{equation*}
Solving these scalar conditions gives $E(\gamma(t))=\e(t^2-2t)/2$. 
According to this, we define
\begin{equation*}
 \gamma^{(\e)}(t):=\gamma(t)-\frac\e2(t^2-2t),
 \qquad \varphi^{(\e)}:=\varphi+\frac\e2.
\end{equation*}
Then $E(\gamma^{(\e)}(t))=0$, $\gamma^{(\e)}(0)=u_0$ and $\gamma^{(\e)}(1)=\varphi^{(\e)}$. Set
\begin{equation*}
 Q:=\int_0^1\norm{\dot\gamma^{(\e)}(t)}_{2,\gamma^{(\e)}(t)}^2dt
\end{equation*}
and denote the terminal velocity by $v:=\dot\gamma(1)$. 

\begin{prop}\label{prop:energy-bounds}
For $0\leq s\leq1$, put
\begin{equation*}
 \sigma(s):=\norm{S_{u_0}-\bar S-s(\tr_{\omega_{u_0}}\chi-\bar\chi)}_{2,u_0}.
\end{equation*}
Every smooth solution of \eqref{eq:cylinder}--\eqref{eq:endpoint} satisfies
\begin{align}
 \norm{v}_{2,\varphi}^2
 &=Q-2\e\int_0^1(I-J)_{\gamma^{(\e)}(t)}(\varphi^{(\e)})dt
 \leq Q,\label{eq:endpoint-kinetic}\\
 \sqrt Q&\leq2\tau\sup_{0\leq s\leq1}\sigma(s)+2C_1\e.
 \label{eq:Q-bound}
\end{align}
It also satisfies the comparison
\begin{equation}\label{eq:endpoint K-bound}
 \begin{split}
 &\tau(M_s(\varphi^{(\e)})-M_s(u_0))+\frac Q2\\
 &\qquad+\e\int_0^1(I-J)_{\gamma^{(\e)}(t)}(\varphi^{(\e)})dt
 \leq\e(I-J)_{u_0}(\varphi^{(\e)}).
 \end{split}
\end{equation}
The constant $C_1$ depends only on the metric comparison in Proposition~\ref{prop:dp}. 
As a consequence, for fixed $u_0,\tau,\chi$, the quantities $d_2(u_0,\varphi^{(\e)})$, $M_s(\varphi^{(\e)})$, $H(\varphi)$ and $\sup_X\varphi$ have upper bounds independent of $s\in[0,1]$ and $0<\e<\min\{1,\tau \Lambda\}$.
\end{prop}
\begin{proof}
Differentiation and spatial integration by parts yield the identity 
\begin{align*}
 \frac d{dt}\frac1{2V}\int_X\dot\gamma^2\omega_\gamma^n
 &=\frac1V\int_X\dot\gamma
              (\ddot\gamma-\abs{\partial\dot\gamma}_{\omega_\gamma}^2)\omega_\gamma^n\\
 &=\e\frac d{dt}\frac1V\int_X\gamma(t)\omega_\varphi^n.
\end{align*}
The second equality uses \eqref{eq:cylinder}. 
Since $\int_X \dot\gamma^{(\e)} \omega_{\gamma^{(\e)}(t)}^n=0$, $\tfrac12\norm{\dot\gamma^{(\e)}}^2-\e\frac1V\int_X\gamma^{(\e)}\omega_\varphi^n$ is as well constant. Evaluating at $t=1$ and using \eqref{eq:I-J} gives
\begin{equation*}
 \frac12\norm{\dot\gamma^{(\e)}(t)}^2
 =\frac12\norm{v}^2
      +\e(I-J)_{\gamma^{(\e)}(t)}(\varphi^{(\e)}).
\end{equation*}
Integration proves \eqref{eq:endpoint-kinetic}.

Let us apply Lemma \ref{lem:energy-comparison} to the comparison path
$\beta(t):=u_0+\e(t^2-2t)/2$. 
Contribution of $\e$-term to kinetic energy is the same for both paths and equals $\tfrac12\int_0^1\e^2(t-1)^2dt=\e^2/6$. Thus
\begin{equation*}
 \tau(M_s(\varphi^{(\e)})-M_s(u_0))+Q/2
 \leq\e\int_0^1\frac1V\int_X(u_0-\gamma^{(\e)}(t))\omega_\varphi^n dt.
\end{equation*}
Since $E(u_0)=E(\gamma^{(\e)}(t))
=E(\varphi^{(\e)})=0$ for every $t$, we may apply
\eqref{eq:I-J} to both terms in the last integral.
The last double integral therefore equals 
$(I-J)_{u_0}(\varphi^{(\e)})-\int_0^1(I-J)_{\gamma^{(\e)}(t)}(\varphi^{(\e)})dt$.
This proves \eqref{eq:endpoint K-bound}.

Propositions~\ref{prop:mabuchi-support} and \ref{prop:dp}, followed by the length--energy inequality, imply
\begin{align*}
 &M_s(\varphi^{(\e)})-M_s(u_0)\geq-\sigma(s)d_2(u_0,\varphi^{(\e)})
                                          \geq-\sigma(s)\sqrt Q,\\
 &(I-J)_{u_0}(\varphi^{(\e)})
 \leq\frac1V\int_X\abs{u_0-\varphi^{(\e)}}\omega_\varphi^n
 \leq C_1d_1(u_0,\varphi^{(\e)})\leq C_1\sqrt Q.
\end{align*}
Substitution into \eqref{eq:endpoint K-bound} gives
$Q/2\leq(\tau\sigma(s)+C_1\e)\sqrt Q$, proving \eqref{eq:Q-bound}.
Dropping the nonnegative terms also gives
\begin{equation*}
 M_s(\varphi^{(\e)})\leq M_s(u_0)+\frac{C_1\e}{\tau}\sqrt Q.
\end{equation*}
On the resulting $d_1$ ball in $E^{-1}(0)$, the nonentropy terms of \eqref{eq:chen-tian}, $J_\chi$, and the spatial supremum are bounded by the energy comparisons of \cite[Lemma 4.4]{CC18}. The identity
\begin{equation*}
 H(\varphi^{(\e)})=M_s(\varphi^{(\e)})-sJ_\chi(\varphi^{(\e)})
 -\bar S E(\varphi^{(\e)})+E_{\Ric(\omega_0)}(\varphi^{(\e)})
\end{equation*}
then bounds entropy. Constant addition preserves entropy so we obtain the entropy bound for $\varphi=\varphi^{(\e)}-\e/2$. 
\end{proof}

\subsection{The explicit solution at $s=1$}
\begin{prop}\label{prop:initial-solution}
At $s=1$, the unique smooth solution of \eqref{eq:cylinder}--\eqref{eq:endpoint} is
\begin{equation}\label{eq:initial-solution}
 \varphi=u_0-\e/2,\qquad \gamma(t)=u_0+\e(t^2-2t)/2.
\end{equation}
\end{prop}
\begin{proof}
Direct computation shows that this pair actually solves \eqref{eq:cylinder}--\eqref{eq:endpoint}. 
It remains to prove the uniqueness. 
Recall the explicit choice of $\Lambda$ in \eqref{eq:twist-choice}. The entropy expression relative to $\omega_{u_0}$ gives
\begin{equation}\label{eq:initial-entropy}
 M_1(\varphi)-M_1(u_0)=H_{u_0}(\varphi)+\Lambda(I-J)_{u_0}(\varphi). 
\end{equation}
Indeed the differentials of the two terms on the right are represented by
$-S_\varphi+\tr_{\omega_\varphi}\Ric(\omega_{u_0})$ and
$\Lambda(\tr_{\omega_\varphi}\omega_{u_0}-n)$. Their sum is the differential of $M+J_\chi$, since $\bar\chi=\bar S+\Lambda n$. Both sides vanish at $u_0$.
Now we substitute \eqref{eq:initial-entropy} into \eqref{eq:endpoint K-bound}. It yields
\begin{equation*}
 \tau H_{u_0}(\varphi^{(\e)})+(\tau \Lambda-\e)(I-J)_{u_0}(\varphi^{(\e)})
 +\frac Q2+\e\int_0^1(I-J)_{\gamma^{(\e)}(t)}(\varphi^{(\e)})dt\leq0.
\end{equation*}
Relative entropy is nonnegative, as is $I-J$, and $\e<\tau \Lambda$. Thus $Q=0$. The path $\gamma^{(\e)}$, whose energy vanishes at every time, is constant. Restoring the energy normalization gives \eqref{eq:initial-solution}. 
\end{proof}

\section{A priori $C^0$-estimates for the endpoint metric}\label{sec:C0}
The first estimate adapts the argument in \cite[Section 5]{CC21I}; see also \cite[Section 3.1]{CC18}. 
In the following lemma, we assume only an upper bound for $f$,
rather than a bound for $\abs{f}$. 
We repeat the argument for the convenience of  readers. 

\begin{lem}\label{lem:CC}
Suppose $\varphi,F,f$ are smooth, $\sup_X\varphi=0$, and
\begin{equation}\label{eq:CC}
      \left\{ 
      \begin{aligned}
  &\omega_\varphi^n=e^F\omega_0^n, \\
&\Delta_{\varphi}F=-f+\tr_{\omega_\varphi}\theta.
      \end{aligned}
      \right.
\end{equation}
Here $\theta$ is a smooth closed real $(1,1)$-form. Assume that entropy in \eqref{eq:chen-tian} is bounded as $H(\varphi)\leq \widehat H$, and that $f\leq M$, $-K\omega_0\leq\theta\leq K\omega_0$. For every $q>1$ there is a bound for $\frac1V\int_Xe^{qF}\omega_0^n$ depending only on $q,\widehat H,M,K$ and the reference geometry. The same data give a bound for $\norm{\varphi}_{C^0}$. 
\end{lem}
\begin{proof}
Fix uniform exponential-integrability constants $\alpha,C_\alpha>0$, so that
\begin{equation}\label{eq:alpha}
 \frac1V\int_Xe^{-\alpha v}\omega_0^n\leq C_\alpha
\end{equation}
for all $v\in\PSH(X,\omega_0)$ with $\sup_Xv=0$. 
This is the integrability used in \cite[Section 5]{CC21I}. The entropy bound and $(-F)e^F\leq e^{-1}$ on $F<0$ imply
\begin{equation*}
 1\leq A:=\frac1V\int_Xe^F\sqrt{1+F^2}\omega_0^n
 \leq1+\widehat H+2/e.
\end{equation*}
Let $\psi$ be the smooth solution of the Monge--Amp\`ere equation 
\begin{equation*}
 \omega_\psi^n=A^{-1}e^F\sqrt{1+F^2}\omega_0^n, 
\end{equation*}
normalized to be $\sup_X\psi=0$.
It exists by famous Calabi--Yau theorem.  
In what follows we take  
\begin{equation*}
 0<\sigma\leq\min\{\frac{1}{2}, \frac{\alpha}{2q}\},\qquad
 \lambda=2(K+2),\qquad \delta=\frac{\alpha}{4n\lambda}. 
\end{equation*}
We will bound $F+\sigma\psi-\lambda\varphi$ from above. Set $U:=\exp(\delta(F+\sigma\psi-\lambda\varphi))$. The arithmetic--geometric mean inequality gives
\begin{equation*}
 \Delta_\varphi\psi\geq nA^{-\frac{1}{n}}(1+F^2)^{\frac{1}{2n}}-\tr_{\omega_\varphi}\omega_0.
\end{equation*}
Using \eqref{eq:CC}, $f\leq M$, and $\Delta_{\varphi}\varphi=n-\tr_{\omega_\varphi}\omega_0$, we obtain
\begin{equation*}
 \delta^{-1}\Delta_\varphi\log U\geq P(F)+c_0\tr_{\omega_\varphi}\omega_0, 
\end{equation*}
where $c_0:=\lambda-K-\sigma>0$ and 
\begin{equation*}
 P(F):=\sigma nA^{-\frac{1}{n}}(1+F^2)^{\frac{1}{2n}} -(M+\lambda n).
\end{equation*}
Take a coordinate ball of fixed small radius $r_0$ centered at a maximum point of $U$. Choose a small constant $\e_0\in(0,1/2)$ and a smooth cutoff $\rho$ on this coordinate ball such that $1-\e_0\leq\rho\leq1$ with $\rho=1$ at its center, $\rho=1-\e_0$ near its boundary, and $\abs{\nabla \rho}\leq C\e_0/r_0$, $\abs{\nabla^2 \rho} \leq C\e_0/r_0^2$. Put $W:=\rho U$. 
The product rule and the identity
$U^{-1}\Delta_\varphi U
=\Delta_\varphi\log U+\abs{\partial\log U}_{\omega_\varphi}^2$
give
\begin{align*}
 U^{-1}\Delta_\varphi W
 &=
 \rho\Delta_\varphi\log U
 +\rho\abs{\partial\log U}_{\omega_\varphi}^2
 +2\operatorname{Re}
   \langle\partial\rho,\partial\log U\rangle_{\omega_\varphi}
 +\Delta_\varphi\rho\\
 &=
 \rho\Delta_\varphi\log U
 +\rho\abs{\partial\log U+\rho^{-1}\partial\rho}
       _{\omega_\varphi}^2
 -\rho^{-1}\abs{\partial\rho}_{\omega_\varphi}^2
 +\Delta_\varphi\rho\\
 &\geq
 \rho\Delta_\varphi\log U
 -\rho^{-1}\abs{\partial\rho}_{\omega_\varphi}^2 
  +\Delta_\varphi\rho. 
\end{align*}
We therefore obtain 

\begin{equation*}
 U^{-1}\Delta_{\varphi}W\geq\delta\rho P(F)+
 \left[\delta\rho c_0-Cr_0^{-2}
 \left(\e_0+\frac{\e_0^2}{1-\e_0}\right)\right]\tr_{\omega_\varphi}\omega_0.
\end{equation*}
For a fixed sufficiently small $\e_0>0$ the coefficient of $\tr_{\omega_\varphi}\omega_0$ is nonnegative. Define $P_-(F):=\max\{-P(F),0\}$. Thus $\Delta_{\varphi}W\geq-\delta U P_-(F)$. Moreover, $P_-(F)\leq M+\lambda n$ and $P(F)<0$ forces
\begin{equation*}
 \abs{F}\leq (1+\widehat H+2/e)\left(\frac{M+\lambda n}{\sigma n}\right)^n=:R.
\end{equation*}
Consequently, $P_-(F)$ vanishes outside
$\{x\in X:\abs{F(x)}\leq R\}$.

If one writes $\Delta_\varphi=A^{ij}\partial_i\partial_j$ in real coordinates of dimension $2n$, \eqref{eq:CC} implies $(\det_{\R}A)^{-1}dV\leq Ce^{2F}\omega_0^n$. Now the determinant form of the real Alexandrov--Bakelman--Pucci estimate \cite[Theorem 9.1]{GT01} gives
\begin{equation*}
 \sup W\leq\sup_{\partial B_{r_0}}W+
 C\left(\frac1V\int_{\abs{F}\leq R}U^{2n}e^{2F}\omega_0^n\right)^{\frac{1}{2n}}.
\end{equation*}
Here the bounded coefficient $\delta(M+\lambda n)$ has been absorbed into $C$. On this integration set, $\psi\leq0$ and $2n\delta\lambda=\alpha/2$ give
\begin{equation*}
 U^{2n}e^{2F}
 \leq e^{(2n\delta+2)R}e^{-\alpha\varphi/2}.
\end{equation*}
Its integral is bounded by \eqref{eq:alpha}, independently of $\sup U$. Since $W=U$ at the maximum point and $W\leq(1-\e_0)\sup U$ at the boundary, the ABP inequality absorbs the boundary term and bounds $\sup U$. We conclude $F+\sigma\psi-\lambda \varphi\leq C$.

As $\varphi\leq0$, this implies $F\leq C-\sigma\psi$. Applying \eqref{eq:alpha} to $\psi$ and using $q\sigma\leq\alpha/2$ gives $V^{-1}\int_Xe^{qF}\omega_0^n\leq e^{qC}C_\alpha$. The compact Monge--Amp\`ere $L^p$ estimate of Ko\l odziej, with $p=2$, bounds the normalized potential $\varphi$. This is the same $L^p$-to-$C^0$ theorem used in \cite[Section 5]{CC21I}. 
\end{proof}

We apply Lemma \ref{lem:CC} to \eqref{eq:cylinder}--\eqref{eq:endpoint}. At the terminal boundary let the inward coordinate $\rho:=1-t$. The explicit subsolution
\begin{equation*}
 \beta(x,t):=\varphi(x)+(1-t)\inf_X(u_0-\varphi)-\frac\e2t(1-t)
\end{equation*}
has spatial metric $\omega_\varphi$, vanishing spatial derivative of its time velocity, and $\ddot\beta=\e$. Thus its geodesic curvature  equals the prescribed density. The boundary values satisfy $\beta(1)=\varphi$ and $\beta(0)\leq u_0$, so comparison principle yields $\beta\leq\gamma$. Convexity also gives $\gamma(t)\leq(1-t)u_0+t\varphi$. Subtracting the common terminal $\phi$ and dividing by $\rho$ we obtain 
\begin{equation*}
 \inf_X(u_0-\varphi)-\frac\e2(1-\rho)
 \leq\frac{\gamma(1-\rho)-\varphi}{\rho}
 \leq u_0-\varphi.
\end{equation*}
Taking $\rho \to 0$ proves
\begin{equation}\label{eq:velocity-barrier}
 \varphi-u_0\leq v:=\dot\gamma(1)\leq\sup_X(\varphi-u_0)+\e/2.
\end{equation}
Set $F:=\log(\omega_\varphi^n/\omega_0^n)$ and
\begin{equation*}
 f:=\bar S-s\bar\chi+v/\tau,\qquad
 \theta:=\Ric(\omega_0)-s\chi.
\end{equation*}
The identity $S_\varphi=-\Delta_{\varphi}F+\tr_{\omega_\varphi}\Ric(\omega_0)$ translates \eqref{eq:endpoint} into \eqref{eq:CC}. Proposition \ref{prop:energy-bounds} bounds $H(\varphi)$.  The same proposition also bounds $\sup \varphi$ so the upper velocity barrier 
in \eqref{eq:velocity-barrier} gives an upper bound for $f$. On the other hand, the lower velocity barrier gives $f\geq \varphi/\tau-C$. 
Since $E(\varphi)=-\e/2$ and
$E(\varphi)\leq\sup_X\varphi$, we have
\begin{equation*}
 -\e/2\leq\sup_X\varphi\leq C.
\end{equation*}
Subtracting $\sup_X\varphi$ changes neither the metric
nor $F$, and changes the constant in
$f\geq\varphi/\tau-C$ by a bounded amount.
Lemma \ref{lem:CC} therefore bounds
$\norm{\varphi-\sup_X\varphi}_{C^0}$, and hence $\norm{\varphi}_{C^0}$. 
Equation~\eqref{eq:velocity-barrier} then bounds $v$,
and hence bounds $f$ from both sides. 

\begin{prop}[Application of Chen--Cheng {\cite[Theorem 5.3]{CC21I}}]\label{prop:CC-higher}
For smooth solutions of \eqref{eq:cylinder}--\eqref{eq:endpoint}, the following estimates hold uniformly in $s\in[0,1]$ and in the allowed range of $\e$. 
There exist $0<\kappa\leq1$ and 
a constant $C_p$ for every $1<p<\infty$ such that 
\begin{equation*}
 \kappa\omega_0\leq\omega_\varphi\leq\kappa^{-1}\omega_0, \qquad 
 \norm{\varphi}_{W^{4,p}}+\norm{F}_{W^{2,p}}\leq C_p.
\end{equation*}
The constants are determined by $\omega_0, u_0, \chi$, and $\tau$. In particular $\varphi$ is bounded in $C^{3,\alpha}$ and $F$ in $C^{1,\alpha}$ for each H\"older exponent $\alpha<1$. \end{prop}
\begin{proof} Proposition~\ref{prop:energy-bounds}, Lemma \ref{lem:CC}, and \eqref{eq:velocity-barrier} verify the entropy bound and the two-sided bounds of $f$. Applying the cited theorem, and then using Sobolev embedding with $p$ arbitrarily large, we conclude the proof. 
\end{proof}

\section{Boundary mixed derivative estimates}\label{sec:cylinder-estimates}

The analysis of non-degenerate complex Monge--Amp\`ere equation on a cylindrical domain appears in 
\cite{Chen00}. The corresponding Dirichlet theory on K\"ahler
manifolds with boundary was also discussed by Phong--Sturm
\cite{PS09}, Guan \cite{Guan09}, and
B{\l}ocki \cite{Blocki12}. Guan--Li
\cite[Theorem 1.1]{GL10} treat the more general Hermitian
setting. 
The gradient maximum principle appears in 
Hanani \cite[Lemma 6]{Hanani96}. 
Our calculation follows 
B{\l}ocki's for compact case 
\cite[proof of Theorem 1, equations (1.14)--(1.16)]{Blocki09};
see also the ``with boundary'' formulation in
\cite[Proposition 3.1]{GL10}.
For mixed boundary derivatives, we use Guan's subsolution
barrier \cite[Section 2]{Guan98}, following the boundary
estimates of Caffarelli--Kohn--Nirenberg--Spruck
\cite{CKNS85} on a strictly pseudoconvex domain. It is used on K\"ahler manifolds in
\cite{Guan09,Blocki12}. 

We repeat their estimates in order to record dependence on the
endpoint $\varphi$ and the prescribed density $\e \omega_\varphi^n/\omega_0^n=:e^H$  varying in the coupled problem.  
As a result, these estimates use
only first derivatives of the logarithmic density $H$,
whereas the interior trace estimate uses second derivatives
of $H$. Remark \ref{rem:bootstrap-order} explains the order
in which they are applied to the coupled system.

\subsection{The Dirichlet problem and its subsolution}
Let
\begin{equation*}
 \mathcal{Z}:=X\times\bigl(\{w\in\C:0\leq\RePart w\leq1\}/\sqrt{-1}\mathbb Z\bigr),
 \qquad w=t+\sqrt{-1}\vartheta. 
\end{equation*}
The quotient by $\sqrt{-1}\mathbb Z$ makes the imaginary direction periodic and the cylinder $\mathcal{Z}$ compact. Let $p_1 \colon \mathcal Z\to X$ be the projection. We consider positive $(1,1)$-forms  
 \begin{align*}
 &\Omega_0:=p_1^*\omega_0+dd^c(t^2/2), \ \text{and} \\
 &\Omega:=\Omega_0+dd^c\Psi, 
 \end{align*}
where $\Psi(x,t):=\gamma(t)(x)-t^2/2$. 
Write $G_0,G$ for the corresponding local coefficient matrices, with the common factor $\frac{\sqrt{-1}}{2\pi}$ taken outside. Thus $\Omega=\frac{\sqrt{-1}}{2\pi}G_{i\bar j}dz^i\wedge d\bar z^j$. 
As a notation, full indices $i,j,k$ range from $1$ to $n+1$; spatial indices $\alpha,\beta$ range from $1$ to $n$, and $z^{n+1}=:w$. Equation \eqref{eq:cylinder} is exactly
\begin{equation}\label{eq:dirichlet}
 \begin{cases}
 \Omega^{n+1}=e^H\Omega_0^{n+1}, \\
 \Psi\vert_{t=0}=u_0,\ \ \text{and} \ \ \Psi\vert_{t=1}=\varphi-\frac{1}{2}.
 \end{cases}
\end{equation}
Here we put 
\begin{equation*}
H(x,t):=\log\e+F(x)= \log \e +\log(\omega_\varphi^n/\omega_0^n). 
\end{equation*}
Indeed $G_{w\bar w}=\ddot\gamma/4$, $G_{\alpha\bar w}=\partial_\alpha\dot\gamma/2$, and the Schur determinant formula gives
\begin{equation*}
 \frac{\det G}{\det G_0}
 =\left(\ddot\gamma-\abs{\partial\dot\gamma}_{\omega_\gamma}^2\right)
      \frac{\omega_\gamma^n}{\omega_0^n}=\e e^F=e^H.
\end{equation*}

Suppose $u_0$ and $\varphi$ have bounded $C^3$-norms and their spatial metrics have a positive lower bound. For sufficiently large $K>0$, define
\begin{equation}\label{eq:strict-subsolution}
 \underline\gamma(t):=(1-t)u_0+t\varphi-Kt(1-t),
 \qquad \underline\Psi:=\underline\gamma-t^2/2.
\end{equation}
Its spatial block is $g_t:=(1-t)g_{u_0}+tg_\varphi$, mixed block is $\partial_\alpha(\varphi-u_0)/2$, and time block is $K/2$. Its Schur complement is
\begin{equation*}
 \frac K2-\frac14g_t^{\alpha\bar\beta}
       \partial_\alpha(\varphi-u_0)\partial_{\bar\beta}(\varphi-u_0).
\end{equation*}
Choosing $K$ large makes this uniformly positive and makes $\det\underline G$ exceed $e^H\det G_0$. Thus $\underline\Omega:=\Omega_0+dd^c\underline\Psi\geq\delta_0\Omega_0$ for some $\delta_0>0$. Comparison and convexity in $t$ give
\begin{equation*}
 \underline\gamma(t)\leq\gamma(t)\leq(1-t)u_0+t\varphi.
\end{equation*}
These inequalities control the potential and its boundary normal derivative. Spatial derivatives at the boundary are prescribed, so they control the full boundary gradient.

\subsection{The gradient estimate}

\begin{prop}[Gradient estimate:
Hanani {\cite[Lemma 6]{Hanani96}} and
B{\l}ocki {\cite[Theorem 1]{Blocki09}}]
\label{prop:gradient}
For a smooth solution of \eqref{eq:dirichlet}, the gradient $\norm{\nabla\Psi}_{C^0(\mathcal Z,G_0)}$ is bounded by a constant depending on $\norm{\Psi}_{C^0}$, the boundary gradient, $\sup H$, and $\norm{dH}_{C^0}$. No second derivative of $H$ is required.
\end{prop}
\begin{proof}

We set $U:=\abs{\nabla\Psi}_{G_0}^2$ and  apply the maximum principle to  $\log U-p(\Psi)$, with weight $p(x):=ax-x^2/2$. At a maximum with $U\geq1$, take background normal holomorphic coordinates which simultaneously diagonalizes $G$, with eigenvalues $\lambda_i$. The critical-point equation and $\Psi_{i\bar j}=(\lambda_i-1)\delta_{ij}$ give
\begin{equation*}
 \sum_j\Psi_{ji}\Psi_{\bar j}
   =(U p'-\lambda_i+1)\Psi_i.
\end{equation*}
Cauchy--Schwarz therefore yields
\begin{equation}\label{eq:ij derivative}
 \frac{1}{U}\sum_{i,j}\lambda_i^{-1}\abs{\Psi_{ji}}^2
 \geq\sum_i\lambda_i^{-1}\abs{\Psi_i}^2
 \left((p')^2-\frac{2p'(\lambda_i-1)}U
                  +\frac{(\lambda_i-1)^2}{U^2}\right).
\end{equation}
Below the $(p')^2$ term will cancel the negative gradient-square term in $\Delta_G\log U$. 
Differentiating
$\log\det G=H+\log\det G_0$ at the maximum point gives
\begin{equation*}
 \sum_i\lambda_i^{-1}\Psi_{i\bar i j}=H_j,
\end{equation*}
since the first derivatives of $G_0$ vanish in the chosen
normal coordinates. Differentiating
$U=\sum_{k,l}(G_0)^{k\bar l}\Psi_k\Psi_{\bar l}$
twice therefore gives
\begin{align*}
 \Delta_G U
 =
 \sum_{i,j}\lambda_i^{-1}
 \left(\abs{\Psi_{ji}}^2+\abs{\Psi_{j\bar i}}^2\right)
 +2\RePart\sum_k H_k\Psi_{\bar k}
 + \sum_{i,k,l}\lambda_i^{-1}
 \partial_i\partial_{\bar i}(G_0)^{k\bar l}
 \Psi_k\Psi_{\bar l}. 
\end{align*}
The second term is estimated as $2\RePart\sum_k H_k\Psi_{\bar k} 
\geq-2\abs{\nabla H}_{G_0}\sqrt U$ 
and the last term is bounded by $C U\tr_G G_0$. 
Dividing by $U$ and discarding the nonnegative mixed
second-derivative terms yields
\begin{equation*}
 \frac{\Delta_G U}{U}
 \geq
 \frac1U\sum_{i,j}\lambda_i^{-1}\abs{\Psi_{ji}}^2
 -C\tr_G G_0
 -\frac{2\abs{\nabla H}_{G_0}}{\sqrt U}.
\end{equation*}
At the maximum point, $U_i=Up'\Psi_i$, and hence
\begin{equation*}
 \Delta_G\log U
 =\frac{\Delta_G U}{U}
 -(p')^2\sum_i\lambda_i^{-1}\abs{\Psi_i}^2.
\end{equation*}
The preceding \eqref{eq:ij derivative} cancels this
last term and gives
\begin{align*}
 \Delta_G\log U
 \geq{}&
 -\frac{2p'}U
 \sum_i\frac{\lambda_i-1}{\lambda_i}\abs{\Psi_i}^2
 +\frac1{U^2}
 \sum_i\frac{(\lambda_i-1)^2}{\lambda_i}\abs{\Psi_i}^2\\
 &-C\tr_G G_0
 -\frac{2\abs{\nabla H}_{G_0}}{\sqrt U}.
\end{align*}
We take $a$ large enough that $p'>0$ on the range of
$\Psi$. Since
\begin{equation*}
 \sum_i\frac{\lambda_i-1}{\lambda_i}\abs{\Psi_i}^2
 =U-\sum_i\lambda_i^{-1}\abs{\Psi_i}^2
 \leq U,
\end{equation*}
the first term on the right is at least $-2p'$.
Discarding the nonnegative square term, we obtain
\begin{equation*}
 \Delta_G\log U
 \geq-2p'-C\tr_G G_0
 -\frac{2\abs{\nabla H}_{G_0}}{\sqrt U}.
\end{equation*}
Finally,
\begin{equation*}
 \Delta_G\Psi=n+1-\tr_G G_0, \ \ \text{so} \ \ 
 \Delta_G p(\Psi)
 =p'(n+1-\tr_G G_0)+p''\abs{\nabla\Psi}_G^2.
\end{equation*}
Subtracting the latter identity gives
\begin{equation}\label{eq:gradient-Laplacian}
 \Delta_G(\log U-p(\Psi))
 \geq(p'-C)\tr_G G_0
 -p''\abs{\nabla\Psi}_G^2-(n+3)p'
 -\frac{2\abs{\nabla H}_{G_0}}{\sqrt U}.
\end{equation} 
Choose $a$ further large so that $p'>C+1$ on the bounded range of $\Psi$. Since $p''=-1$, the maximum principle in \eqref{eq:gradient-Laplacian} bounds both the inverse trace and $\abs{\nabla\Psi}_G^2$ at this point. 
Let us compare $\abs{\nabla\Psi}_G^2$ with $U=\abs{\nabla\Psi}_{G_0}^2$. 
At the maximum point, choose $A\geq1$ such that
\begin{equation*}
 \sum_{i=1}^{n+1}\lambda_i^{-1}\leq A,
 \qquad
 \abs{\nabla\Psi}_G^2\leq A. 
\end{equation*}
The first inequality gives $\lambda_i\geq A^{-1}$ for every $i$.
Since the Monge--Amp\`ere equation gives
$\prod_{i=1}^{n+1}\lambda_i=e^H$, we obtain
\begin{equation*}
 \lambda_i
 =\frac{e^H}{\prod_{j\ne i}\lambda_j}
 \leq A^n\sup_{\mathcal{Z}}e^H.
\end{equation*}
In the chosen coordinates, it follows that
\begin{equation*}
 U=\sum_i\abs{\Psi_i}^2
 \leq
 \left(\max_i\lambda_i\right)
 \sum_i\lambda_i^{-1}\abs{\Psi_i}^2
 \leq A^{n+1}\sup_{\mathcal{Z}}e^H.
\end{equation*}
Finally, 
the bounded range of $p(\Psi)$ and the boundary gradient estimate give a global gradient bound. 
\end{proof}

\subsection{Mixed derivatives at  boundaries}

\begin{prop}[
{\cite[Section 2]{Guan98}}]
\label{prop:mixed}
Assume that \eqref{eq:dirichlet} has the strict subsolution with a uniform positive lower bound and a $C^3$ bound (see \eqref{eq:strict-subsolution}). A $C^1$ bound for $H$ and the gradient estimate in Proposition~\ref{prop:gradient} give
\begin{equation}\label{eq:mixed-bound}
 \sup_{\partial \mathcal{Z}}\abs{\nabla_X\partial_t\Psi}_{G_0}\leq C.
\end{equation}
The constant does not require a bound for $\nabla^2H$. 
\end{prop}
\begin{proof}
Fix a spatial coordinate ball $B_R(p) \subset X$ and the inward coordinate $\rho=t$ at the initial boundary $t=0$. 
We work on the tube $B_R(p)\times[0,1/5]\times(\R/\mathbb Z)$. Write
\begin{equation*}
 \mathscr L:=G^{i\bar j}\partial_i\partial_{\bar j},\qquad
 T:=\sum_{i=1}^{n+1} G^{i\bar i},\qquad h:=\Psi-\underline\Psi\geq0.
\end{equation*}
The strict subsolution gives $\underline G\geq2\delta_0 I$ and
$\mathscr Lh=n+1-\tr(G^{-1}\underline G)$. For positive constants $a,b$, set $h_b:=h+a\rho-2b\rho^2$. Then, since $\mathscr L(\rho^2)=G^{w\bar w}/2$ in the cylinder coordinate $w=\rho+i\vartheta$,
\begin{equation*}
 \mathscr Lh_b\leq n+1-2\delta_0T-bG^{w\bar w}.
\end{equation*}
Put $K_b:=\diag(\delta_0, \dots, \delta_0, b)$. The arithmetic--geometric mean inequality gives
\begin{equation*}
 \tr(G^{-1}K_b)\geq (n+1)\left(\frac{\det K_b}{\det G}\right)^{\frac{1}{n+1}}.
\end{equation*}
An upper bound for $\det G$ permits a fixed choice of $b$ for which the right-hand side is at least $n+2$. With $a=2b/5+5$ this gives
\begin{equation}\label{eq:collar}
 \mathscr Lh_b\leq-1-\delta_0T,\qquad
 h_b\geq0,\qquad h_b\vert_{\rho=\frac{1}{5}}\geq1,\qquad h_b\vert_{\rho=0}=0.
\end{equation}

For a real spatial derivative $D$ in each coordinate, differentiation of the determinant equation yields 
\begin{equation*}
 \mathscr L(Dh)
 =D\log\det G-
 G^{i\bar j}\big(D(G_0)_{i\bar j}
                   +\partial_i\partial_{\bar j}D\underline\Psi\big).
\end{equation*}
We derive $\abs{\mathscr L(Dh)}\leq C_0+C_1T$ from it, using $\norm{dH}_{C^0}$ 
and subsolution $C^3$ bounds. The global gradient estimate gives $\abs{Dh}\leq K_D$. The spatial function $\abs{x-p}^2$ satisfies $0\leq\mathscr L\abs{x-p}^2\leq C_1T$. Choose $cR^2\geq K_D$, and then take $A$ large enough such that
\begin{equation*}
 B_\pm:=Ah_b+c\abs{x-p}^2\pm Dh
\end{equation*}
are nonnegative on the lateral boundary
$\partial B_R(p)\times[0,1/5]\times(\R/\mathbb{Z})$ and on
$B_R(p)\times\{1/5\}\times(\R/\mathbb{Z})$, and satisfy
$\mathscr L B_\pm<0$. 
On $B_R(p)\times\{0\}\times(\R/\mathbb{Z})$, we have $Dh=0$, and hence
$B_\pm=c\abs{x-p}^2\geq0$.
These sets form the boundary of the tube.
The minimum principle therefore gives $B_\pm\geq0$
throughout the tube. 
Taking inward derivatives at $(p,0)$ yields
\begin{equation*}
 \abs{\partial_\rho Dh(p,0)}
 \leq A((\partial_\rho h)(p,0)+a)
 \leq A(\norm{\Psi}_{C^1}+\norm{\underline\Psi}_{C^1}+a).
\end{equation*}
Adding back the controlled derivative of the subsolution bounds the initial mixed derivatives.

For the terminal boundary use the holomorphic coordinate $w'=1-w$, whose inward real coordinate is $\sigma=1-t$. The same $h$ vanishes there and is nonnegative. Define $h_b=h+a\sigma-2b\sigma^2$ on the terminal tube. It satisfies \eqref{eq:collar} with $\sigma$ in place of $\rho$, and $Ah_b+c\abs{x-p}^2\pm Dh$ have the same boundary signs and negative $\mathscr L$. The inward derivative estimate therefore bounds $\partial_{\sigma}Dh$ at $\sigma=0$, which is $-\partial_tDh$ at $t=1$. 
\end{proof}

\subsection{Solvability of the  Dirichlet problem}
\begin{prop}[Application of Guan--Li {\cite[Theorem 1.1]{GL10}}]\label{prop:dirichlet}
For smooth $u_0,\varphi\in\cH$ and $\e>0$, \eqref{eq:dirichlet} has a unique smooth solution. It is invariant under translations of $\vartheta$. 
For fixed $\e>0$, every integer $k\geq4$, and
$0<\alpha<1$, the terminal velocity $v=\dot{\gamma}(1)$ depends smoothly on $\varphi \in C^{k,\alpha}$ satisfying
$\omega_\varphi>0$, with values in $C^{k-1,\alpha}(X)$, 
as described in \eqref{eq:velocity-map}.
\end{prop}
\begin{proof}
We first apply 
\cite[Theorem 1.1]{GL10}. The product cylinder
$(\mathcal Z,\Omega_0)$ is a compact K\"ahler manifold
with smooth boundary. For fixed smooth endpoints and
$\e>0$, the density $e^H$ is smooth and strictly positive.
The function $\underline\Psi$ in
\eqref{eq:strict-subsolution} has the prescribed boundary
values and satisfies
\begin{equation*}
 \Omega_0+dd^c\underline\Psi>0,\qquad
 (\Omega_0+dd^c\underline\Psi)^{n+1}
 \geq e^H\Omega_0^{n+1}.
\end{equation*}
Thus the hypotheses of the theorem are satisfied, and
a smooth solution exists. Uniqueness follows from comparison principle. 

We next record the second-order estimates. 
Here we need a control of $\nabla^2H$. 
For the interior trace $\tr_{G_0}G$, the Aubin--Yau
calculation gives
\begin{equation}\label{eq:cylinder-trace}
 \Delta_G(\log \tr_{G_0}G-A\Psi)
 \geq\frac{\Delta_{G_0}H-C}{\tr_{G_0}G}
         +(A-C)\tr_G G_0-(n+1)A.
\end{equation}
See \cite{Yau78}. 
If one takes $A$ large and observes 
$\tr_{G_0}G \geq (n+1)e^{\frac{1}{n+1}H}$, this bounds $\tr_{G}G_0$ at an interior maximum of $\log \tr_{G_0}G-A\Psi$. 
Similarly to the proof of Proposition \ref{prop:gradient}, 
setting $\lambda_i$ as the eigenvalues relative to $G_0$, the inequalities $\sum_i\lambda_i^{-1}\leq A$ and $\prod_i\lambda_i\leq e^H$ imply
\begin{equation*}
 \lambda_i=\frac{\prod_j\lambda_j}{\prod_{j\ne i}\lambda_j}\leq A^ne^H,
 \qquad \tr_{G_0}G\leq (n+1)A^ne^H.
\end{equation*}
Boundary second derivatives are controlled by \eqref{eq:mixed-bound} and the Schur complement formula: 
\begin{equation*}
 G_{w\bar w}
 =\frac{\det G}{\det G_{\rm sp}}
   +\sum_{\alpha,\beta=1}^nG_{\alpha\bar w}G_{w\bar\beta}(G_{\rm sp}^{-1})^{\bar\beta\alpha}, 
\end{equation*}
where $G_{\rm sp}:=(G_{\alpha\bar\beta})_{1\leq\alpha,\beta\leq n}$ is the spatial block. 
Therefore we obtain $\tr_{G_0}G$ bound up to the boundary. 
Combining it with the determinant lower bound finally gives a lower bound for each $\lambda_i$.  

For $k\geq4$ and $0<\alpha<1$, the implicit function
theorem gives local smooth dependence of the solution
on $\varphi$ in the open set
\begin{equation*}
 \{\varphi\in C^{k,\alpha}(X):\omega_\varphi>0\}.
\end{equation*}
The prescribed logarithmic density belongs to
$C^{k-2,\alpha}$, and the solution belongs to
$C^{k,\alpha}$. Thus the map 
\begin{equation}\label{eq:velocity-map}
 v:\varphi\longmapsto\dot\gamma_\varphi(1)
\end{equation}
is locally $C^1$ with values in $C^{k-1,\alpha}(X)$.
It is uniformly bounded on families satisfying
$\norm{\varphi}_{C^{k,\alpha}}\leq B$ and
$\omega_\varphi\geq\kappa\omega_0$, for fixed
$B$ and $\kappa>0$.
Note that the bounds may depend on $\e$. 
\end{proof}

\begin{prop}\label{prop:all-bounds}
Fix $u_0,\tau,\e>0$ as in Section~\ref{sec:coupled}. For every integer $k\geq4$ and $0<\alpha<1$, there are $\kappa_0>0$ and $B_0<\infty$, independent of $s\in[0,1]$ and of the smooth solution of \eqref{eq:cylinder}--\eqref{eq:endpoint}, such that
\begin{equation}\label{eq:all-bounds}
 \omega_\varphi\geq\kappa_0\omega_0,\qquad
 \norm{\varphi}_{C^{k,\alpha}}\leq B_0.
\end{equation}
\end{prop}
\begin{proof}
Proposition~\ref{prop:CC-higher} gives $\varphi\in C^{3,\alpha}$ and $F\in C^{1,\alpha}$ with uniform bounds. The strict subsolution \eqref{eq:strict-subsolution} and Propositions~\ref{prop:gradient} and \ref{prop:mixed} therefore bound $\nabla v$ of $v=\dot{\gamma}(1)$. 
We next apply this bound to the prescribed scalar curvature
equation. 
The identity 
$S_\varphi=\tr_{\omega_\varphi}\Ric(\omega_0)-\Delta_\varphi F$
and the boundary equation \eqref{eq:endpoint} give
\begin{equation*}
 \Delta_\varphi F
 =-\bar S+s\bar\chi-\frac{v}{\tau}
  +\tr_{\omega_\varphi}\bigl(\Ric(\omega_0)-s\chi\bigr).
\end{equation*}
The gradient bound makes $v$ uniformly Lipschitz, while the
already established $C^{3,\alpha}$ bound for $\varphi$ and
uniform metric equivalence of Proposition \ref{prop:CC-higher} control the trace term in $C^{1,\alpha}$. 
The Poisson equation in \eqref{eq:CC} now has a uniformly H\"older right-hand side and uniformly elliptic $C^{1,\alpha}$ coefficients. Schauder estimates give $F\in C^{2,\alpha}$. Differentiating the endpoint Monge--Amp\`ere equation gives $\varphi\in C^{4,\alpha}$. 
The logarithmic density $H=\log\e+F$ in \eqref{eq:dirichlet} now has its required second-derivative bound, so \eqref{eq:cylinder-trace} applies. 

Higher regularity of $\gamma$ controls the terminal velocity $v$.
The prescribed scalar curvature and endpoint Monge--Amp\`ere equations 
then improve the regularity of $F$ and $\varphi$, which in turn 
improves the density and boundary data for the equation \eqref{eq:dirichlet}. 
Iterating these estimates gives \eqref{eq:all-bounds}
at every order. 
\end{proof}

\begin{rem}\label{rem:bootstrap-order}
The estimates are applied in the following order:
\begin{align*}
 (\varphi,F)\in C^{3,\alpha}\times C^{1,\alpha}
 &\quad\text{by Proposition~\ref{prop:CC-higher}},\\
 \norm{\nabla\Psi}_{L^\infty}\leq C
 &\quad\text{by Proposition~\ref{prop:gradient}},\\
 \norm{\nabla v}_{L^\infty}\leq C
 &\quad\text{by Proposition~\ref{prop:mixed}},\\
 (\varphi, F)\in C^{4,\alpha}\times C^{2,\alpha}
 &\quad\text{by the Schauder estimates in Proposition~\ref{prop:all-bounds}},\\
 G\asymp G_0
 &\quad\text{by \eqref{eq:cylinder-trace} and its boundary estimates},\\
 (\varphi,\gamma)\in C^\infty
 &\quad\text{by the bootstrap in Proposition~\ref{prop:all-bounds}}.
\end{align*}

The interior trace estimate uses second derivatives of $H$. 
Since $H=\log\e+F$, applying that estimate before
obtaining $F\in C^{2,\alpha}$ would be circular.

The Dirichlet estimates are uniform for families of data
with uniform bounds on all the quantities on which their
constants depend. In the coupled problem, however, the
endpoint $\varphi$ is unknown and the density
$e^H=\e\omega_\varphi^n/\omega_0^n$ depends on it.
Fixed-endpoint solvability alone therefore does not
supply the required uniform bounds. 

Finally, $\e$ and $\tau$ are fixed during continuation
in $s$. The cylinder ellipticity constants may
deteriorate as $\e\to 0$. 
In
Section~\ref{sec:recovery}, 
we instead use the energy comparison  
when taking the limit. 
\end{rem}

\section{Leray--Schauder degree argument}
For a bounded open set $\mathscr D$ in a real Banach space and a continuous compact (not necessarily linear) map $T:\overline{\mathscr D}\to\mathscr B$, with $x\ne T(x)$ on $\partial\mathscr D$, the Leray--Schauder degree $\deg_{\rm LS}(I-T,\mathscr D,0)$ is an integer. It is defined by finite-dimensional approximation and Brouwer degree; see \cite[Chapter 2, Sections 8.3--8.6]{Deimling85}. We use three properties. A nonzero degree implies a fixed point. If an open subset $U\subset\mathscr D$ contains all fixed points, excision gives
\begin{equation*}
 \deg_{\rm LS}(I-T,\mathscr D,0)=\deg_{\rm LS}(I-T,U,0).
\end{equation*}
Finally, if $T_s$ is continuous jointly in $(s,x)$, the set
$\bigcup_{s\in[0,1]}T_s(\overline{\mathscr D})$ has compact closure, and there are no boundary fixed points, then the degree is independent of $s$. These properties and the local linearization formula are recalled in \cite[Theorem 1.1 and the subsequent discussion]{LS22}.

For fourth-order equations related to scalar curvature, Chau--Weinkove \cite[Section 2, proof of Theorem 1.1]{CW15} use a compact fixed-point homotopy for a second boundary value problem, following Trudinger--Wang \cite[Lemma 5.2]{TW05}. Their boundary problem differs from our equation. We specify the compact map and its initial degree below.

Let us define  
\begin{equation*}
 G_s(\varphi):=\tau[S_\varphi-\bar S-s(\tr_{\omega_\varphi}\chi-\bar\chi)]-v(\varphi),
 \qquad v(\varphi)=\dot\gamma_\varphi(1),
\end{equation*}
where $v$ is the map in \eqref{eq:velocity-map}. Uniqueness of the Dirichlet problem gives $\gamma_{\varphi+c}(t)=\gamma_\varphi(t)+tc$, so $G_s(\varphi+c)=G_s(\varphi)-c$. In particular the equation determines the additive constant. 
Fix $k\geq6$ and $0<\alpha<1$. The  little H\"older space is
\begin{equation*}
 h^{k,\alpha}(X):=\overline{C^\infty(X,\R)}^{C^{k,\alpha}(X)}.
\end{equation*}
Choose $0<\kappa<\kappa_0$ and $B>B_0$ in Proposition~\ref{prop:all-bounds}, and set
\begin{equation*}
 \mathscr D:=\{\varphi\in h^{k,\alpha}(X):\omega_\varphi>\kappa\omega_0,
                       \ \norm{\varphi}_{C^{k,\alpha}}<B\}.
\end{equation*}
 Relative to a fixed smooth connection $\nabla$, write
\begin{equation*}
 -G_s(\varphi)=A_\varphi\nabla^4\varphi+R_s(\varphi),
 \qquad \sigma_4(A_\varphi\nabla^4)(x,\xi)
       =\tau(g_\varphi^{i\bar j}\xi_i\bar\xi_j)^2.
\end{equation*}
Here $A_\varphi\nabla^4\varphi=A_\varphi^{ijkl}\nabla_i\nabla_j\nabla_k\nabla_l \varphi$ denotes the contraction of tensor $A_\varphi$ with $\nabla^4 \varphi$, and $\sigma_4$ denotes the principal symbol. Expanding $\Delta_\varphi\log\det g_\varphi$ leaves only products of derivatives of order at most three after this fourth-order term is removed. 
The term $\tr_{\omega_\varphi}\chi$ depends on derivatives
of $\varphi$ of order at most two, while
\eqref{eq:velocity-map} maps $C^{k,\alpha}$ potentials
to $C^{k-1,\alpha}$ functions.
Consequently, $R_s(\varphi)$ is bounded in
$C^{k-3,\alpha}$ whenever
$\norm{\varphi}_{C^{k,\alpha}}\leq B$ and
$\omega_\varphi\geq\kappa\omega_0$,
with a bound uniform in $s\in[0,1]$. 
The operator $A_\varphi\nabla^4$ is  uniformly strongly elliptic.  
The tensor $A_\varphi$ is bounded
in $C^2$ on $\overline{\mathscr D}$: the lower bound $\omega_\varphi\geq\kappa\omega_0$ and the uniform $C^4$ bound for $\varphi$ control the inverse metric and its first two derivatives. Fix a smooth volume form $dV$, say $\omega_0^n$, 
and use it for all Sobolev norms and pairings in this paragraph. 
Two integrations by parts give
\begin{equation*}
 \langle A_\varphi\nabla^4f,f\rangle
 =\int_X A_\varphi^{ijkl}(\nabla_k\nabla_lf)
   (\nabla_i\nabla_jf)dV+\mathcal R_\varphi(f),
\end{equation*}
where derivatives of the coefficients and of $dV$, and commutators
of the fixed connection, contribute only lower-order terms. It holds 
\begin{equation*}
 \abs{\mathcal R_\varphi(f)}
 \leq C(\norm{f}_{H^2}\norm{f}_{H^1}+\norm{f}_{H^1}^2) 
\end{equation*}
uniformly in
$\varphi\in\overline{\mathscr D}$. 
Since the principal symbol was the square of that of the Laplacian $\Delta_\varphi$, it is natural to  compare 
the integrand with $\tau\Delta_\varphi^2$. 
Further integrations 
by parts give 
\begin{equation*}
 \langle A_\varphi\nabla^4f,f\rangle
 =\tau\norm{\Delta_\varphi f}_{L^2(dV)}^2+\widetilde{\mathcal R}_\varphi(f), 
\end{equation*}
where the residual term still satisfies 
\begin{equation*}
 \abs{\widetilde{\mathcal R}_\varphi(f)}
 \leq C(\norm{f}_{H^2}\norm{f}_{H^1}+\norm{f}_{H^1}^2).
\end{equation*}
Note that $\Delta_\varphi$ is no longer self-adjoint with respect to the fixed volume form $dV$. 
Interpolation 
$ \norm{f}_{H^1}^2\leq\eta\norm{f}_{H^2}^2+C_\eta\norm{f}_{L^2}^2$ 
combined with 
Young's inequality bounds the remainder by 
$\eta'\norm{f}_{H^2}^2+C_{\eta'}\norm{f}_{L^2}^2$ for any $\eta'>0$. 
From the uniform second-order elliptic estimate
\begin{equation*}
 \norm{f}_{H^2}^2\leq C(\norm{\Delta_\varphi f}_{L^2}^2+\norm{f}_{L^2}^2)
\end{equation*}
we take $\eta'<\tau/(2C)$ such that 
\begin{equation}\label{eq:garding}
 \RePart\langle A_\varphi\nabla^4f,f\rangle
 \geq c_0\norm{f}_{H^2}^2-C_0\norm{f}_{L^2}^2
\end{equation}
holds for all $\varphi\in\overline{\mathscr D}$. 
Choose $K>C_0$ and put
$P_\varphi:=A_\varphi\nabla^4+K$. Lax--Milgram theorem and elliptic Schauder regularity give inverses
\begin{equation}\label{inverse}
 P_\varphi^{-1}\colon h^{k-4,\alpha}\longrightarrow h^{k,\alpha},\qquad
 P_\varphi^{-1}\colon h^{k-3,\alpha}\longrightarrow h^{k+1,\alpha}, 
\end{equation}
with operator norms uniformly bounded in $\varphi \in \overline{\mathscr D}$.  
The coefficient regularity for the second assertion is available because $A_\varphi\in C^{k-2,\alpha}$. 
Since $P_{\varphi_j}^{-1}-P_{\varphi}^{-1}=
P_{\varphi_j}^{-1}(P_\varphi-P_{\varphi_j})P_\varphi^{-1}$, 
the inverses depend
continuously on $\varphi$ in the indicated operator norms. 
Using the second inverse in \eqref{inverse}, we now define 
\begin{equation*}
 T_s(\varphi):=P_\varphi^{-1}(K\varphi-R_s(\varphi)).
\end{equation*}
A fixed-point of $T_s$ is precisely $G_s(\varphi)=0$. 
The preceding uniform estimates give
\begin{equation*}
 \sup_{s\in[0,1],\ \varphi\in\overline{\mathscr D}}
       \norm{T_s(\varphi)}_{C^{k+1,\alpha}}<\infty. 
\end{equation*}
The compact embedding 
$h^{k+1,\alpha}\Subset h^{k,\alpha}$ follows from Arzel\`a--Ascoli theorem and H\"older interpolation $
[f]_{C^{0,\alpha}}
\leq C\norm{f}_{C^0}^{1-\alpha}
\norm{f}_{C^1}^{\alpha}$. 
Consequently $T_s$ is a uniformly compact homotopy, in the precise sense stated above. 
The fixed-point property enables us to 
iterate the regularity gain so that every fixed point is smooth. 
As we took $\kappa <\kappa_0$ and $B>B_0$, Proposition~\ref{prop:all-bounds} then excludes fixed points on $\partial\mathscr D$.

Let us calculate the degree at $s=1$. In the rest of this calculation, $\varphi=u_0-\e/2$ is the solution obtained in \eqref{eq:initial-solution}. For an endpoint variation $h \in C^\infty(X; \R)$, define the variation of $\gamma$ as 
\begin{equation*}
 f(t):=D_h[\varphi \mapsto \gamma_\varphi(t)]=\left.\frac d{da}\right|_{a=0}\gamma_{\varphi+ah}(t).
\end{equation*}
Here and below $D$ denotes the Fr\'echet derivative of a map between the indicated Banach spaces. 
Linearizing \eqref{eq:cylinder} at \eqref{eq:initial-solution} gives
\begin{equation*}
 \ddot f=\e\Delta_{u_0}(h-f),\qquad f(0)=0,\quad f(1)=h.
\end{equation*}
For eigenvector $-\Delta_{u_0}h=\lambda_jh$ and $k_j:=\sqrt{\e\lambda_j}$, its solution and the derivative of the boundary velocity map are
\begin{equation*}
 f(t)=\left(1-\frac{\sinh(k_j(1-t))}{\sinh k_j}\right)h,
 \qquad D_hv(\varphi)=\dot f(1)=\frac{k_j}{\sinh k_j}h.
\end{equation*}
For $\lambda_j=0$, these formulas mean $f(t)=th$ and $D_hv(\varphi)=h$. 
To compute the linearization of 
\eqref{eq:endpoint} at \eqref{eq:initial-solution}, we first calculate the Hessian of $M_1$ at $u_0$ using
\eqref{eq:initial-entropy}. 
The first variation of the density is $\Delta_{u_0}h$, so relative entropy has second variation $V^{-1}\int_X(\Delta_{u_0}h)^2\omega_{u_0}^n$. The variation formula for $I-J$ gives its second variation $-V^{-1}\int_Xh\Delta_{u_0}h\omega_{u_0}^n$. Thus
\begin{equation*}
 D^2_{h, h}M_1\big\vert_{u_0}=\frac1V\int_Xh(\Delta_{u_0}^2h-\Lambda\Delta_{u_0}h)\omega_{u_0}^n.
\end{equation*}
The eigenvalues of $-DG_1(\varphi)$ are therefore
\begin{equation}\label{eq:initial-spectrum}
 \tau(\lambda_j^2+\Lambda\lambda_j)+\frac{k_j}{\sinh k_j}>0.
\end{equation}
The constant mode has eigenvalue one. The elliptic operator $\tau(\Delta_{u_0}^2-\Lambda\Delta_{u_0})+I$ is an isomorphism; the difference from $-DG_1(\varphi)$ is compact from $h^{k,\alpha}$ to $h^{k-4,\alpha}$, since $Dv$ loses only one derivative. Hence $DG_1(\varphi)$ is Fredholm of index zero. Formula \eqref{eq:initial-spectrum} proves its kernel is zero, so index zero implies surjectivity and hence invertibility.

The identity $P_\psi(\psi-T_s(\psi))=-G_s(\psi)$, differentiated at a fixed point $\psi$, gives
\begin{equation*}
 I-DT_s(\psi)=-P_\psi^{-1}DG_s(\psi).
\end{equation*}
Here we take the first inverse in \eqref{inverse}. 
The term containing $DP_\psi$ vanishes because $\psi-T_s(\psi)=0$. In particular $I-DT_1(\varphi)$ is invertible. The local linearization theorem for Leray--Schauder degree therefore gives
$\deg_{\rm LS}(I-T_1,U,0)\in\{-1,1\}$ on a sufficiently small neighborhood $U$ of $\varphi$; see \cite[Theorem 1.1 and its linear degree formula]{LS22}. Uniqueness at $s=1$, excision, and homotopy invariance now give
\begin{equation*}
 \deg_{\rm LS}(I-T_0,\mathscr D,0)
 =\deg_{\rm LS}(I-T_1,\mathscr D,0)
 =\deg_{\rm LS}(I-T_1,U,0)\ne0.
\end{equation*}
Thus the coupled system has a solution at $s=0$.

\begin{prop}\label{prop:coupled-existence}
For $u_0\in\cH$ with $E(u_0)=0$, $\tau>0$, and $0<\e<\min\{1,\tau \Lambda\}$, system \eqref{eq:cylinder}--\eqref{eq:endpoint} at $s=0$ has a smooth solution. 
\end{prop}
\begin{proof}
The uniformly compact homotopy $T_s$ has nonzero degree and no boundary fixed points. Its fixed points are precisely the smooth solutions, as shown above. 
\end{proof}

\section{Recovery of the infimum of the Calabi functional}\label{sec:recovery}
\begin{thm}\label{prop:smooth-recovery}
On every compact K\"ahler manifold,
\begin{equation*}
 \inf_{\varphi\in\cH}\Cal(\omega_\varphi)
 =\inf_{u\in\cE^2}\abs{\partial M}(u).
\end{equation*}
\end{thm}
\begin{proof}
Let us fix $u_0,\tau$ and choose solutions $(\varphi_{\tau,\e},\gamma)$ of \eqref{eq:cylinder}-\eqref{eq:endpoint} from Proposition~\ref{prop:coupled-existence}. 
We first observe that normalization $E=0$ is compatible with the K-resolvent. Indeed, every smooth path satisfies
\begin{equation*}
 \norm{\dot\gamma}_{2,\gamma}^2
 =\norm{\dot\gamma-\tfrac d{dt}E(\gamma)}_{2,\gamma}^2
       +\abs{\tfrac d{dt}E(\gamma)}^2.
\end{equation*}
Replacing a path $\gamma(t)$ by
$\gamma(t)-E(\gamma(t))$ makes its Monge--Amp\`ere
energy zero at every time.
Among real-valued paths joining $E(u)$ to $E(v)$,
the affine path minimizes the integral of the squared
derivative. Taking the infimum of the above
decomposition over paths joining $u$ to $v$ we therefore obtain  
\begin{equation*}
 d_2(u,v)^2=d_2(u-E(u),v-E(v))^2+\abs{E(u)-E(v)}^2.
\end{equation*}
Completion extends this formula to $\cE^2$. Since $M(u+c)=M(u)$ for every constant $c\in\R$, it follows that the resolvent $p_\tau$ of $u_0$ has $E(p_\tau)=0$. 

Next fix $\delta>0$ and a smooth $w$ with $E(w)=0$, which will be taken as an approximant to $p_\tau$.
 Choose a smooth path $\alpha\colon [0,1]\to\cH$ with $E(\alpha(t))=0$ from $u_0$ to $w$ such that
\begin{equation*}
 \int_0^1\norm{\dot \alpha(t)}_{2,\alpha(t)}^2dt\leq d_2(u_0,w)^2+\delta.
\end{equation*}
Existence follows from the definition of the path-length distance, reparametrization 
\begin{equation*}
t\mapsto \left(\int_0^1\sqrt{\norm{\dot{\alpha}(t)}_{\alpha(t)}^2+\eta^2}dt\right)^{-1}\int_0^t \sqrt{\norm{\dot{\alpha}(\tau)}_{\alpha(\tau)}^2+\eta^2}d\tau 
\end{equation*}
for a small $\eta>0$, 
 and subtraction of $E$, which decreases kinetic energy. Now we apply Lemma \ref{lem:energy-comparison} to 
$\beta(t):=\alpha(t)+\e(t^2-2t)/2$.
As in the proof of Proposition~\ref{prop:energy-bounds}, the contributions of $p(t)=\e(t^2-2t)/2$ cancels so that 
\begin{equation}\label{eq:proximal-comparison}
 \begin{split}
 &\tau M(\varphi^{(\e)})+\frac Q2
       +\e\int_0^1(I-J)_{\gamma^{(\e)}(t)}(\varphi^{(\e)})dt
 \\ &\quad\leq\tau M(w)+\frac{d_2(u_0,w)^2+\delta}{2}
       +\e\int_0^1(I-J)_{\alpha(t)}(\varphi^{(\e)})dt.
 \end{split}
\end{equation}
For the fixed path $\alpha$, we bound the last integral uniformly
in $\e$. Since $E(\alpha(t))=E(\varphi^{(\e)})=0$,
\eqref{eq:I-J} gives
\begin{equation*}
 (I-J)_{\alpha(t)}(\varphi^{(\e)})
 =\frac1V\int_X(\alpha(t)-\varphi^{(\e)})\omega_\varphi^n.
\end{equation*}
 Proposition~\ref{prop:dp} with $p=1$, the inequality $d_1\leq d_2$,
and the length--energy inequality imply
\begin{align*}
 \frac1V\int_X\abs{\varphi^{(\e)}-u_0}\omega_\varphi^n
 \leq C_1d_1(u_0,\varphi^{(\e)})
 \leq C_1d_2(u_0,\varphi^{(\e)})
 \leq C_1\sqrt Q.
\end{align*}
Consequently,
\begin{equation*}
 0\leq (I-J)_{\alpha(t)}(\varphi^{(\e)})
 \leq \norm{\alpha(t)}_{C^0}+\norm{u_0}_{C^0}+C_1\sqrt Q.
\end{equation*}
The smooth path $\alpha$ is fixed, and \eqref{eq:Q-bound} bounds
$\sqrt Q$ independently of $\e$ for fixed $\tau$. Thus there is
a constant $C_\alpha$, independent of $\e$, such that
\begin{equation*}
 0\leq \e\int_0^1(I-J)_{\alpha(t)}(\varphi^{(\e)})dt
 \leq  C_\alpha\e. 
\end{equation*}
We define, in the normalization of Lemma \ref{lem:resolvent},
\begin{equation*}
 F_\tau(u):=M(u)+\frac{d_2(u,u_0)^2}{2\tau},\quad
 m:=F_\tau(p_\tau),\quad A_\e:=M(\varphi^{(\e)})+\frac Q{2\tau}.
\end{equation*}
To estimate $A_\e$, 
divide \eqref{eq:proximal-comparison} by $\tau$ and drop the nonnegative integral on the left. 
First let $\e\to 0$, then $\delta\to0$, finally by Proposition~\ref{prop:entropy-recovery} approximate $p_\tau$ by smooth potentials $w_j$ satisfying $E(w_j)=0$. These limits give $\limsup_{\e\to0}A_\e\leq m$. Conversely, minimality and the length--energy inequality give
\begin{equation*}
 m\leq F_\tau(\varphi^{(\e)})\leq A_\e.
\end{equation*}
Both quantities therefore tend to $m$. Strong convexity gives
\begin{equation*}
 d_2(\varphi^{(\e)},p_\tau)^2
 \leq2\tau(F_\tau(\varphi^{(\e)})-m)\longrightarrow0,
\end{equation*}
and subtracting the two defining energies gives
\begin{equation*}
 0\leq Q-d_2(\varphi^{(\e)},u_0)^2
 =2\tau(A_\e-F_\tau(\varphi^{(\e)}))
 \leq2\tau(A_\e-m)\longrightarrow0.
\end{equation*}
Thus $Q\to d_2(p_\tau,u_0)^2$. At $s=0$, \eqref{eq:endpoint} and \eqref{eq:endpoint-kinetic} imply
\begin{equation*}
 \Cal(\omega_\varphi)=\frac{\norm{v}_{2,\varphi}}\tau
 \leq\frac{\sqrt Q}\tau,
 \qquad
 \limsup_{\e\to0}\Cal(\omega_{\varphi_{\tau,\e}})
 \leq\frac{d_2(p_\tau,u_0)}\tau.
\end{equation*}
Letting $\tau\to\infty$ and using Lemma \ref{lem:resolvent} proves
$\inf_{\cH}\Cal\leq\inf_{\cE^2}\abs{\partial M}$. 
\end{proof}

\section{Algebraic approximation with control of the second moment}\label{sec:algebraic}
\subsection{Non-Archimedean metrics and approximation}
From now on we fix an ample line bundle $L$. 
Let $X^{\NA}$ be the Berkovich analytification of $X$. A divisorial valuation of the function field $\C(X)$ naturally defines a point 
of this space. In particular, $v_{\rm triv}$ denotes the trivial valuation. 
Divisorial valuations only form part of $X^\NA$ but 
play a dominant role below. 
Following the works of Boucksom--Jonsson, 
we first explain how a given test configuration defines a metric on $X^{\NA}$. 
We use the conventions of \cite[Sections 4.2--4.7]{BJ25}. 
In this aspect, it is natural to allow $\cL$ appearing in the definition of a test configuration
not necessarily ample but relatively semiample ($L$ is assumed to be ample throughout).
This is because if one takes a common equivariant birational model 
\begin{equation*}
p \colon \overline{\cY} \longrightarrow \overline{\cX}, \qquad \rho \colon \overline{\cY} \longrightarrow X \times \mathbb{P}^1, 
\end{equation*} 
the pullback polarization $p^*\overline{\cL}$ on
$\overline{\cY}$ defines a semiample test configuration.
Taking such a model, 
we may assume $(\cX, \cL)$ itself dominates $X \times \C$. 
Write $L_{\cX}:=\rho^*p_1^*L$ for the pullback of $L$ from the product. 
Using the fixed equivariant identification over $\C^*$, the line bundle can be uniquely written as $\cL=L_{\cX}+D$ with a $\Q$-divisor $D$ such that $\supp D \subseteq \cX_0$. 

Let us take a such model $p\colon \cY\to\cX$. 
For each irreducible component $F$ of $\cY_0$, we define a valuation $v_F:=\ord_F(\cY_0)^{-1}\ord_F\vert_{\C(X)}$. 
A real-valued function $f=\varphi_{(\cX, \cL)}$ on $X^{\NA}$ is defined as 
\begin{equation*}
 f(v_F):=\frac{\ord_F(p^*D)}{\ord_F(\cY_0)}. 
\end{equation*}
This is independent of the choice of models.
The above definition is justified because any divisorial valuation with rational coefficient is equivalent to some $v_F$ obtained in this way \cite[Theorem 6.27]{BHJ17}
and the above $f$ is 
determined by divisorial valuations \cite[Theorem 4.6]{BHJ17}. 
We call the above $f$ {\em a model function}. 
If $(\cX, \cL)$ is a semiample model, 
then $f$ defines a continuous plurisubharmonic function on $X^\NA$ (strictly speaking, relative to $L$). 
We call it {\em semiample model function}. 
As a matter of fact, we define a non-Archimedean plurisubharmonic function to be an upper-semicontinuous function obtained as a decreasing limit of such semiample model functions; see \cite[Section 4.4]{BJ25}. 
Thus the symbol $f$ can be regarded as a non-Archimedean analogue of 
the K\"ahler potentials $\varphi$ and 
we sometimes call a non-Archimedean psh function a metric. 
The trivial model metric defined by the product $X \times \mathbb{P}^1$ corresponds to the reference metric $\omega_0$ in the K\"ahler setting. 

After passing to a dominating smooth model, assume that
$(\cX, \cX_{0,\red})$ is log smooth, namely $\cX_{0, \red}$ is  simple normal crossing. 
Its dual complex is then defined as 
\begin{equation*}
 \Delta_{\cX}
 :=
 \left(
 \coprod_{\emptyset\ne J\subset I}
 \coprod_{Z\subset \cap_{i\in J}E_i}
 \left\{(w_i)_{i\in J}\in\R_{\geq0}^J:
               \sum_{i\in J}b_iw_i=1\right\}
 \right)\bigg/{\sim},
\end{equation*}
where the faces are identified according to inclusions
of strata.
We write $\cX_0 =\sum_i b_i E_i$ and $Z$ runs through irreducible components of $\cap_{i \in J}E_i$. 
Even if $\cL$ is not relatively semiample, model functions are still continuous and piecewise linear on $\Delta_{\cX} \subset X^{\NA}$, 
however, they need not be plurisubharmonic in $X^{\NA}$.

Set  $b_E:=\ord_E(\cX_0)$. 
The Monge--Amp\`ere energy 
and operator for semiample model metrics are defined as 
\begin{align}\label{eq:na-energy-ma}
 &E^{\NA}(\varphi_{(\cX, \cL)}):=\frac{\overline{\cL}^{n+1}}{(n+1)V}, \\ 
 &\MA^{\NA}(\varphi_{(\cX, \cL)})
 :=\frac1V\sum_E b_E(\cL^n\cdot E)\delta_{v_E}. 
\end{align}
We extend the definition of $E^{\NA}$ to arbitrary plurisubharmonic 
functions by taking the above decreasing limit approximation. 
It allows $E^{\NA}(\varphi) =-\infty$, similarly to the K\"ahler setting. 
The class $\cE^{1,\NA}$ consists of those non-Archimedean plurisubharmonic functions with finite energy. 
We can then define the Monge--Amp\`ere operator in $\cE^{1, \NA}$, using {\em any} decreasing sequence of semiample model metrics. 
One can also define the distance $d_1^{\NA}$ on $\cE^{1,\NA}$ by the formula 
\begin{equation*}
 d_1^{\NA}(\varphi,\psi)
 =E^{\NA}(\varphi)+E^{\NA}(\psi)
       -2E^{\NA}(P(\min\{\varphi,\psi\})),
\end{equation*}
where $P(f)$ is a plurisubharmonic envelope. 

The log discrepancy function is first 
defined 
$A_X(c\ord_E)=c(1+\ord_EK_{Y/X})$ for divisorial valuations. 
(This $E$ denotes a prime divisor over $X$.)
It is then defined as the maximal lower-semicontinuous extension to $X^{\NA}$. 
The non-Archimedean entropy is defined as 
\begin{equation*}
 H^{\NA}(\varphi):=\int_{X^{\NA}}A_X\MA^{\NA}(\varphi).
\end{equation*}
Also, the non-Archimedean Ricci energy $R^{\NA}$ is the extension of the intersection number description $V^{-1}(\rho^*p_1^*K_X\overline{\cL}^{n})$ on semiample models. Its mixed-energy formulation extends continuously in the strong topology. 
With these conventions, we obtain the formula 
\begin{equation*}
 M^{\NA}=H^{\NA}+R^{\NA}+\bar S E^{\NA} 
\end{equation*}
as in \cite[Section 7.4]{BHJ17}. 
In particular, it follows that  convergence in $d_1^{\NA}$ and in $H^{\NA}$ implies convergence in $M^{\NA}$.

We next recall the definition and energy-related properties of maximal weak geodesic ray associated with $\varphi \in \cE^{1, \NA}$, as a natural extension of subsection \ref{sec:test configurations}. 
As we have already noted in section \ref{sec:normalization}, {\em weak geodesic} stands for a weak solution of \eqref{eq:geodesic}. 
A finite-energy weak geodesic 
is a weak geodesic in the metric space $(\cE^p,d_p)$. 
Subgeodesics are the subsolutions 
of \eqref{eq:geodesic}. 
A subgeodesic $u_t$ is called ray if it is defined for all $t \in [0, \infty)$ and it is said to have linear growth if $\sup_X u_t \leq C(1+t)$. 
For non-Archimedean $\varphi\in\cE^{1,\NA}$ and a fixed smooth initial value $v$, the associated {\em maximal} weak geodesic ray is characterized by
\begin{equation*}
 \ell_{\varphi, t}
 :={\sup}^*\bigg\{u_t: U\text{ is a subgeodesic ray of linear growth},\
       u_0\leq v,\ U^{\NA}\leq\varphi\bigg\}. 
\end{equation*}
The superscript $*$ denotes upper-semicontinuous regularization. 
It coincides with Berman--Boucksom--Jonsson
\cite[Definition 6.5 and Theorem 6.6]{BBJ21}. 
When $\varphi$ is the model metric of a semiample test configuration, this is equivalent to
the condition explained in \ref{sec:test configurations}: 
the extendability of $U$ to a locally bounded psh weight. 
See \cite[Corollary 5.3]{BJ25}. 
On the other hand, a subgeodesic ray $u_t$ defines the corresponding $\mathbb{S}^1$-invariant
plurisubharmonic potential $U(x, \tau):=u_{-\log{\abs{\tau}}^2}(x)$ on the  punctured product space $X \times \mathbb{B}^*$. 
If it has linear growth, the generic Lelong number $\nu_E(\rho^*U)$
along an irreducible component $E$ of a fixed model is well-defined. 
Given a defining equation $z_E=0$ 
our convention is $\rho^*U=\nu_E(\rho^*U)\log\abs{z_E}^2+O(1)$ near the generic point of $E$. 
We define the non-Archimedean metric 
\begin{equation}\label{eq:ray to NA}
U^{\NA}(v_E): = -\frac{\nu_E(\rho^*U)}{b_E}, 
\end{equation}
which is not necessarily in $\cE^{1, \NA}$. 
If $\varphi=\varphi_{(\cX, \cL)}$ is defined by a semiample model, $U^{\NA}$ of the associated weak geodesic ray $u_t$ coincides with $\varphi$ and is automatically in $\cE^{1, \NA}$. 
The following results are also crucial in recent progress of Yau--Tian--Donaldson conjecture. 
Recall that using convexity we define the radial energy as 
\begin{equation*}
M^{\rad}(\ell):= \lim_{t \to \infty} \frac{M(\ell_t)}{t}. 
\end{equation*}

\begin{prop}[Li {\cite[Theorem 1.2]{Li22}};
Boucksom--Jonsson {\cite[Theorem 8.7]{BJ25}}]
\label{prop:na-slope}
Any weak geodesic ray with finite radial K-energy is maximal.
For every $\varphi\in\cE^{1,\NA}$, its associated
maximal ray $\ell_\varphi$ satisfies
\begin{equation}\label{eq:na-slope}
 M^{\NA}(\varphi)=M^{\rad}(\ell_\varphi).
\end{equation}
The equality is understood in $\mathbb R\cup\{+\infty\}$.
\end{prop}

For rays $\ell,\ell'$ the chordal distance can be similarly defined as 
\begin{equation*}
 d_p^c(\ell,\ell'):=\lim_{t\to\infty}\frac{d_p(\ell_t,\ell'_t)}t.
\end{equation*}
For a constant-speed ray $\ell$ in $\cE^p$, we denote $\norm{\ell}_p:=d_p(\ell_0,\ell_1)$. 

\begin{prop}[Maximal-ray isometry {\cite[Section 5.1]{BJ25}}]\label{prop:ray-isometry}
For a fixed initial potential, the correspondence
$\varphi\mapsto\ell_\varphi$ preserves distance:
\begin{equation*}
 d_1^c(\ell_\varphi,\ell_\psi)
 =d_1^{\NA}(\varphi,\psi).
\end{equation*}

\end{prop}

Given a valuation $v \in X^{\NA}$, 
we denote $G(v)$ by the Gauss extension of $v$, which is defined rational functions in $\C(\cX)=\C(X \times \C)$ as 
\begin{equation*}
 G(v)\left(\sum_{k} a_k\tau^k\right):=\min_k\big\{v(a_k)+k\big\}. 
\end{equation*}
This is normalized to $G(v)(\tau)=1$. 
For a fixed model $\cX$, 
Non-zero values $w_i:=G(v)(E_i)$ for irreducible components define a point $r(v)$ in $\Delta$ and the natural retraction 
\begin{equation*}
r\colon X^{\NA} \longrightarrow \Delta_{\cX}. 
\end{equation*}

\begin{prop}[Boucksom--Jonsson {\cite{BJ22, BJ25}}]\label{prop:retraction} \noindent
      \begin{itemize}
      \item[$(1)$]
Let $\varphi\in\cE^{1,\NA}$ and index the simple normal crossing models by domination: $j\preceq k$ means that model indexed $k$ dominates model $j$. We write $\Delta_j$ for the dual complex and $r_j\colon X^{\NA}\to\Delta_j$ for the natural retraction. There is a unique $\varphi_j\in\cE^{1,\NA}$ with $\sup\varphi_j=\sup\varphi$ satisfying
\begin{equation}\label{eq:retracted-measures}
 \MA^{\NA}(\varphi_j)=(r_j)_*\MA^{\NA}(\varphi).
\end{equation}
Each $\varphi_j$ is continuous, and the net $\varphi_j$ converges to $\varphi$ in $d_1^{\NA}$ and in entropy. If $H^{\NA}(\varphi)<\infty$, these convergences can be realized by a sequence. 
\item[$(2)$]
Let $\Delta$ be the dual complex of a fixed snc model,
with retraction $r$. 
For a 
$\psi\in\cE^{1,\NA}$ satisfying $ \supp\MA^{\NA}(\psi)\subseteq\Delta$, $ \sup_{X^{\NA}}\psi=0$, 
there exists rational piecewise-linear functions $f_i$ on
$\Delta$ such that $P(f_i\circ r)$ converges to $\psi$
uniformly and in entropy.
Each $f_i \circ r$ defines a model function with $D:=\sum_E b_E f_i(r(v_E))E$, 
where we take a dominating model realizing a subdivision on which $f_i$  is affine. 
\end{itemize}
\end{prop}
\begin{proof}
Let $\mathcal M^{1,\NA}$ denote the space of finite-energy probability measures, endowed with the strong topology. 
\cite[Theorem 12.8]{BJ22} states that
\begin{equation*}
 \MA^{\NA}\colon \cE^{1,\NA}/\R\longrightarrow\mathcal M^{1,\NA}
\end{equation*}
is a homeomorphism, 
as an analogue of \cite[Theorem 4.7]{BBGZ13}. 
It guarantees existence of $\varphi_j$. 
Since $(r_j)_*\MA^{\NA}(\varphi)$ is supported on $\Delta_j$, \cite[Theorem 12.12, Example 12.13]{BJ22} makes $\varphi_j$ continuous. 
By \cite[Lemma 4.10]{BJ25}, the measures $(r_j)_*\MA^{\NA}(\varphi)$ belong to $\mathcal M^{1,\NA}$, satisfy $(r_j)_*\MA^{\NA}(\varphi)\longrightarrow\MA^{\NA}(\varphi)$ strongly as a net,  and 
\begin{equation*}
\int A_X(r_j)_*\MA^{\NA}(\varphi)\longrightarrow\int A_X\MA^{\NA}(\varphi) 
\end{equation*}
holds. 
The continuity of the inverse 
of $\MA^{\NA}$
then gives $d_1^{\NA}(\varphi_j,\varphi)\to0$. 
When $H^{\NA}(\varphi)$ is finite, 
one can choose successively a dominating model for which both $d_1^{\NA}$ and the entropy errors are controlled. This gives the required sequence. 

For a fixed dual complex, \cite[Lemmas 7.16 and 7.17]{BJ25} give the envelope approximation and entropy convergence. The proof first approximates $f:=\psi\vert_\Delta$ by rational piecewise-linear functions and then uses
\begin{equation*}
 \sup_{X^{\NA}}\abs{P(f_i\circ r)-P( f\circ r)}\leq \sup_{\Delta}\abs{f_i-f}. 
\end{equation*}
which is \cite[Section 5.3, (5.1)]{BJ22}. 
When the Monge--Amp\`ere measure is supported on $\Delta$, \cite[Proposition 6.24(i)]{BJ25} assures $\psi=P(\psi\vert_{\Delta}\circ r)$ so this gives an approximation of $\psi$. 
\end{proof}

We also use the following consequence of hybrid convergence explained in \cite[secion 5.4]{BJ25}. 
Let $\psi_0,\ldots,\psi_m$ be semiample model metrics, and let
$u_{0,t},\ldots,u_{m,t}$ be compatible smooth
rays. 
For every continuous function $F\colon\mathbb R^m\to\mathbb R$,
\begin{equation*}
 \lim_{t\to\infty}
 \int_X
 F\left(\frac{u_{1, t}}{t},\ldots,
        \frac{u_{m, t}}{t}\right)
 \MA\bigl(u_{0, t}\bigr)
 =
 \int_{X^{\NA}}F(\psi_1,\ldots,\psi_m)
 \MA^{\NA}(\psi_0).
\end{equation*}
In particular, for a smooth ray $u_t$ compatible 
with semiample model metric $\psi$ and any $p\geq1$, we obtain
\begin{equation}\label{eq:higher moment}
 \lim_{t\to\infty}
 \frac{1}{Vt^p}\int_X\abs{u_t}^p\omega_{u_t}^n
 =
 \int_{X^{\NA}}\abs{\psi}^p\MA^{\NA}(\psi). 
 \end{equation} 
 and 
\begin{equation}\label{eq:hybrid sup}
 \lim_{t \to \infty} \frac{\sup_{X} u_t}{t} = \sup_{X^{\NA}} \psi, \qquad
 \lim_{t \to \infty} \frac{\inf_{X} u_t}{t} = \inf_{X^{\NA}} \psi. 
\end{equation}
Indeed, the rescaled potentials extend continuously to the
non-Archimedean fiber metric of the hybrid space, with values
$\psi_1,\ldots,\psi_m$, while $\MA(u_{0, t})$
converges weakly to $\MA^{\NA}(\psi_0)$;
see \cite[Section 5.4, Definition 5.12 and Example 5.13]{BJ25}.
The rescaled potentials are uniformly bounded for $t\geq1$ so only the restriction of $F$
to a fixed compact subset of $\mathbb R^m$ is relevant. 
Together with Lemma \ref{lem:radial-moment}, the formula \eqref{eq:higher moment} allows us to estimate the velocity moments in Proposition \ref{prop:dh-measure} through mixed Monge--Amp\`ere integrals.

\subsection{A mixed-moment bound}
In order to show a uniform higher-moment bounds for the solutions 
to \eqref{eq:retracted-measures}, 
We provide the following. 
Set 
\begin{equation*}
 e_p(\psi):=\int_{X^{\NA}}(-\psi)^p\MA^{\NA}(\psi)
\end{equation*}
for $\psi\in\cE^{1,\NA}$ with
$\sup_{X^{\NA}}\psi=0$, allowing the value $+\infty$. 
Note that in general $\sup_{X^{\NA}} \psi =\psi(v_{\rm triv})$, by \cite[Proposition 4.12 (ii)]{BJ22}. 

\begin{lem}\label{lem:mixed-moment}
Let $\varphi$ be a bounded non-Archimedean plurisubharmonic metric with $-M\leq\varphi\leq0$ and $\sup_{X^{\NA}}\varphi=0$. Let $\psi$ be a continuous such metric with $\sup_{X^{\NA}}\psi=0$. For every integer $p>n$, we obtain 
\begin{equation}\label{eq:mixed-moment}
 \int_{X^{\NA}}(-\psi)^p\MA^{\NA}(\varphi)
 \leq \sum_{k=1}^n\binom nk(p)_kM^k\cdot  
          e_p(\psi)^{\frac{p-k}{p}},  
\end{equation}
where $(p)_k:=p(p-1)\cdots(p-k+1)$.
\end{lem}

\begin{proof}
First suppose both metrics are semiample model metrics. Choose compatible smooth positive rays $u_t,v_t$ and subtract their spatial suprema. This does not change their Monge--Amp\`ere measures or the associated non-Archimedean potentials. The formula \eqref{eq:hybrid sup} gives $-M_t\leq u_t\leq0$, $v_t\leq0$, and $\limsup M_t/t\leq M$. 
For $k+l\leq n$ set
\begin{equation*}
 I(q,k,l)=\frac1V\int_X(-v_t)^q
       \omega_{v_t}^k\wedge\omega_{u_t}^l
                   \wedge\omega_0^{n-k-l}.
\end{equation*} 
Integrating twice by parts gives
\begin{equation}\label{eq:mixed-recursion}
 I(q,k,l)\leq I(q,k,l-1)+qM_t I(q-1,k+1,l-1).
\end{equation}
Indeed, set $T:=\omega_{v_t}^k\wedge\omega_{u_t}^{l-1}
       \wedge\omega_0^{n-k-l}$. 
Since $T$ is closed, integration by parts gives
\begin{align*}
 I(q,k,l)-I(q,k,l-1)
 =\frac1V\int_X(-v_t)^q dd^cu_t\wedge T
 =\frac1V\int_Xu_t dd^c(-v_t)^q\wedge T.
\end{align*}
Using
\begin{equation*}
 dd^c(-v_t)^q
 =-q(-v_t)^{q-1}dd^cv_t
  +q(q-1)(-v_t)^{q-2}dv_t\wedge d^cv_t
\end{equation*}
and $dd^cv_t=\omega_{v_t}-\omega_0$, we obtain
\begin{align*}
 I(q,k,l)-I(q,k,l-1)
 =
 &\frac qV\int_X(-u_t)(-v_t)^{q-1}
                   \omega_{v_t}\wedge T
 -\frac qV\int_X(-u_t)(-v_t)^{q-1}
                   \omega_0\wedge T\\
 &+\frac{q(q-1)}V\int_Xu_t(-v_t)^{q-2}
                   dv_t\wedge d^cv_t\wedge T.
\end{align*}
The last two terms are nonpositive, since $u_t\leq0$
and $T$ is positive. Applying $-u_t\leq M_t$ to the
first term yields
\begin{align*}
 I(q,k,l)-I(q,k,l-1)
 &\leq\frac{qM_t}{V}\int_X(-v_t)^{q-1}
                       \omega_{v_t}\wedge T\\
 &=qM_t I(q-1,k+1,l-1),
\end{align*}
which proves \eqref{eq:mixed-recursion}. 
Iterating over the $n$ factors gives
\begin{equation*}
 I(p,0,n)\leq\sum_{k=0}^n\binom nk(p)_k M_t^k I(p-k,k,0).
\end{equation*}
For $q\geq1$, a further integration by parts shows
$I(q,k,0)\leq I(q,k+1,0)$. Replace each $k\geq1$ term with $I(p-k,n,0)$ and divide by $t^p$. The term $k=0$ tends to zero: normalized $\omega_0$-plurisubharmonic functions have uniformly bounded reference $L^p$ norms by \eqref{eq:alpha}. The other terms converge by \eqref{eq:higher moment}. We obtain
\begin{equation*}
 \int_{X^{\NA}}(-\psi)^p\MA^{\NA}(\varphi)
 \leq\sum_{k=1}^n\binom nk(p)_k M^k e_{p-k}(\psi).
\end{equation*}
H\"older for the probability measure $\MA^{\NA}(\psi)$ bounds
$e_{p-k}(\psi)\leq e_p(\psi)^{(p-k)/p}$, proving \eqref{eq:mixed-moment} for smiample models. 

For a semiample model metric $\varphi$ and continuous
$\psi$, the inequality \eqref{eq:mixed-moment} follows from uniform approximation
of $\psi$ by semiample model metrics. 
For bounded $\varphi$, choose semiample model metrics
$\varphi_i\searrow\varphi$.
Their suprema decrease to zero.
After subtracting these suprema, the approximants lie
in $[-M-\delta, 0]$ for all sufficiently large $i$. 
Their Monge--Amp\`ere measures are unchanged by this
normalization and converge weakly to $\MA^{\NA}(\varphi)$.
Since $(-\psi)^p$ is continuous, we may pass to the limit
in the inequality with $M+\delta$ in place of $M$.
Letting $\delta \to 0$ we complete the proof.  
\end{proof}

\begin{lem}\label{lem:radial-moment}
For a continuous non-Archimedean plurisubharmonic metric $\psi$ with 
$\sup_{X^{\NA}}\psi=0$ and $p\geq1$, the associated maximal weak geodesic ray is contained in $\cE^p$ 
and its radial $d_p$ speed of satisfies
\begin{equation*}
 c_{p}e_p(\psi)\leq\norm{\ell_\psi}_p^p
                \leq C_{p}e_p(\psi).
\end{equation*}
For $p>2$ and $\alpha:=\frac{p-2}{2(p-1)}$, weak geodesic rays in $\cE^p$ also satisfy the interpolation inequality 
\begin{equation}\label{eq:interpolation}
 d_2^c(\ell,\ell')\leq C_{p}
       d_1^c(\ell,\ell')^\alpha d_p^c(\ell,\ell')^{1-\alpha}. 
\end{equation}
\end{lem}
\begin{proof}
First suppose that $\psi$ is a semiample model metric,
and choose a smooth compatible ray $u_t$.
Since $\sup_{X^{\NA}}\psi=0$,
\eqref{eq:hybrid sup} gives 
$t^{-1}\sup_X u_t \to 0$. 
Thus \eqref{eq:alpha} implies $u_t -\sup_X u_t$ has uniform $L^p(\omega_0^n)$-norm. 
Consequently,
\begin{align*}
 \frac1{Vt^p}\int_X\abs{u_t}^p\omega_0^n
 &\leq
 2^{p-1}\left(
 \frac1{Vt^p}\int_X\abs{u_t-\sup_X u_t}^p\omega_0^n
 +\frac{\abs{\sup_X u_t}^p}{t^p}\right)\longrightarrow0.
\end{align*}
On the other hand, \eqref{eq:higher moment} and $\psi\leq0$
give
\begin{equation*}
 \lim_{t\to\infty}
 \frac1{Vt^p}\int_X\abs{u_t}^p\omega_{u_t}^n
 =
 \int_{X^{\NA}}(-\psi)^p\MA^{\NA}(\psi)
 =e_p(\psi).
\end{equation*}
Applying \eqref{eq:dp-comparison} to $u_t$ and $0$, we obtain
\begin{equation*}
 C_p^{-1}\frac{d_p(u_t,0)^p}{t^p}
 \leq
 \frac1{Vt^p}\int_X\abs{u_t}^p
             \left(\omega_{u_t}^n+\omega_0^n\right)
 \leq
 C_p\frac{d_p(u_t,0)^p}{t^p}.
\end{equation*}

Let $\ell_t$ be the maximal weak geodesic ray associated with $\psi$ so that 
$\abs{u_t -\ell_t} \leq C$. 
By \eqref{eq:dp-comparison}, this implies
$d_p(u_t,\ell_t)=O(1)$. 
The triangle inequality gives
\begin{equation*}
 \abs{d_p(u_t,0)
      -d_p(\ell_t,\ell_0)}
 \leq d_p(u_t,\ell_t)+d_p(0,\ell_0)
 =O(1).
\end{equation*}
Since $\ell$ has constant $d_p$ speed,
$t^{-1}d_p(u_t,0) \to \norm{\ell_\psi}_p$. 
The above comparison therefore yields
\begin{equation*}
 C_p^{-1}\norm{\ell_\psi}_p^p
 \leq e_p(\psi)
 \leq C_p\norm{\ell_\psi}_p^p,
\end{equation*}
which proves the first assertion for semiample model metrics.

Now let $\psi$ be continuous.
Choose semiample model metrics $\psi_j$, with 
$\sup_{X^{\NA}}\psi_j=0$, such that
$\psi_j \to \psi$ uniformly. 
Uniform convergence implies weak convergence of the
Monge--Amp\`ere measures, and also 
$(-\psi_j)^p\to(-\psi)^p$ uniformly. Hence the right-hand side of 
\begin{align*}
 \abs{e_p(\psi_j)-e_p(\psi)}
 \leq
 \abs{\int_{X^{\NA}}(-\psi)^p
       \left(\MA^{\NA}(\psi_j)-\MA^{\NA}(\psi)\right)} 
        +\sup_{X^{\NA}} \abs{(-\psi_j)^p-(-\psi)^p }
\end{align*}
converges to zero. 
On the other hand, let us take the associated maximal weak geodesic rays $\ell, \ell_j$ with a common initial point, for $\psi, \psi_j$ respectively.  
Their envelope characterization, together with 
$\psi-\delta_j\leq\psi_j\leq\psi+\delta_j$, gives
\begin{equation*}
 \ell_t-\delta_j t
 \leq\ell_{j,t}
 \leq\ell_t+\delta_jt.
\end{equation*}
Since $\psi$ is continuous, compactness of $X^{\NA}$ implies that these rays have bounded potentials on each finite time interval, and hence take values in $\cE^p$.
Using \eqref{eq:dp-comparison}, we obtain
\begin{equation*}
 d_p(\ell_{j,t},\ell_t)
 \leq (2C_p)^{1/p}\delta_jt.
\end{equation*}
The triangle inequality gives
\begin{equation*}
 \abs{\norm{\ell_j}_p-\norm{\ell}_p}
 \leq (2C_p)^{1/p}\delta_j
 \longrightarrow0.
\end{equation*}
The first assertion for continuous $\psi$ now follows. 

For the second assertion, fix $u,v\in\cE^p$ and set
$\nu:=\MA(u)+\MA(v)$. 
From the choice of $\alpha$ 
H\"older's inequality yields
\begin{equation*}
 \norm{u-v}_{L^2(\nu)}
 \leq
 \norm{u-v}_{L^1(\nu)}^\alpha
 \norm{u-v}_{L^p(\nu)}^{1-\alpha}.
\end{equation*}
Applying \eqref{eq:dp-comparison} three times we obtain 
\begin{align*}
 d_2(u,v)
 &\leq C_2^{1/2}\norm{u-v}_{L^2(\nu)}
 \leq C_2^{1/2}
       \norm{u-v}_{L^1(\nu)}^\alpha
       \norm{u-v}_{L^p(\nu)}^{1-\alpha}\\
 &\leq C_2^{1/2}C_1^\alpha C_p^{(1-\alpha)/p}
       d_1(u,v)^\alpha d_p(u,v)^{1-\alpha}.
\end{align*}
Application of this inequality to $u=\ell_t$ and $ v_t=\ell'_t$ gives 
\begin{equation*}
 \frac{d_2(\ell_t,\ell'_t)}t
 \leq C_p
 \left(\frac{d_1(\ell_t,\ell'_t)}t\right)^\alpha
 \left(\frac{d_p(\ell_t,\ell'_t)}t\right)^{1-\alpha}.
\end{equation*}
It proves \eqref{eq:interpolation}. 
\end{proof}

\subsection{Ample model approximation of a model envelope}
We combine base-ideal approximation of \cite{BJ22} with entropy convergence of \cite[Theorem B]{Tru26}, which is an evolved version of Fujita approximation. 
Uniform convergence of the base-ideal approximants implies $L^2$-convergence. 
This construction over $k((t))$ is already explained in Boucksom\textquotesingle s appendix \cite[Theorem A.6]{Tru26}. 
The argument below is written directly for equivariant models over the trivially valued field $\C$. 

\begin{prop}[Consequence of Boucksom--Jonsson and Trusiani {\cite[Lemma 5.22]{BJ22}, \cite[Theorem B]{Tru26}}]\label{prop:model-approximation}
Let $(\cX, \cL)$ be a smooth simply normal crossing model with $\cL$ a relatively big $\Q$-line bundle.  Let $f=f_{(\cX, \cL)}$ be the model function and set $\psi=P(f)$. There exist normal ample test configurations $(\cX_j,\cL_j)$ with model metrics $\psi_j$, uniform constant $M$ such that
\begin{equation}\label{eq:model-approximation}
 d_1^{\NA}(\psi_j,\psi)\to0,\qquad
 H^{\NA}(\psi_j)\to H^{\NA}(\psi),\qquad
 \abs{\psi_j} \leq M. 
\end{equation}
 The maximal weak geodesic rays $\ell_{\psi_j}, \ell_\psi$ satisfy
 \begin{equation}
 d_2^c(\ell_{\psi_j}, \ell_\psi) \to 0, \qquad 
 M^{\rad}(\ell_{\psi_j}) \to M^{\rad}(\ell_\psi),
 \end{equation} 
 and convergence of the centered norms 
 \begin{equation}
 \norm{\ell_{\psi_j}-E^{\rad}(\ell_{\psi_j})}_2 \to \norm{\ell_\psi-E^{\rad}(\ell_\psi)}_2. 
 \end{equation}
\end{prop}
\begin{proof}
We repeat the arguments of cited papers 
in order to obtain additional convergence results. 
Write $\cL=L_{\cX}+D$ with $\supp D \subseteq \cX_0$. Choose a sufficiently large $C \in \mathbb{Q}$ so that $-C\cX_0\leq D\leq C\cX_0$. For sufficiently divisible $k$, pushforward to $\mP^1$ gives inclusions
\begin{equation*}
 H^0(X,kL)\otimes\cO_{\mathbb{P}^1}(-kC)
 \subseteq \pi_*(k\overline{\cL})
 \subseteq  H^0(X,kL)\otimes\cO_{\mathbb{P}^1}(kC).
\end{equation*}
The second inclusion implies that 
the middle sheaf is torsion-free on the smooth curve, hence it is locally free. 
Let us write the equivariant splitting as 
$\pi_*(k\overline{\cL}) \simeq \oplus_{j}\cO_{\mathbb{P}^1}(\lambda_{j}^{(k)})$. 
The above inclusions imply $-kC\leq \lambda_j^{(k)}\leq kC$. Indeed, the lower bound follows from generic surjectivity and that $\Hom(\cO_{\mathbb{P}^1}(-kC), \cO_{\mathbb{P}^1}(\lambda))\neq 0$ only if $kC+\lambda \geq 0$. 
Let us temporarily replace $\cL$ with $\cL +C\cX_0$. 
The corresponding weights $\lambda_j^{(k)}+kC$ are nonnegative, so $\pi_*(k\overline{\cL})$ is globally generated for every sufficiently divisible $k$. 
In particular, the absolute and relative base ideals agree. 
Take $C$ larger if necessary, so that $\overline{\cL}$ is big and its augmented base locus $\mathbb{B}_+(\overline{\cL})$ is vertical.  
Note that rescaling $\cL \mapsto \cL+C\cX_0$ only adds a constant to $f,\psi$
so it does not affect to the approximation. 

After an equivariant modification supported on the central
fiber, we may also assume that 
$(\cX, \cL)$ is a smooth simply normal crossing model and 
$
 \overline{\cL}=\overline{\cA}_{\cX}+F,
$
where $\overline{\cA}_{\cX}$ is ample and $F$ is an
effective vertical $\Q$-divisor. 
This preparation does not change $f$ or $P(f)$. See the preliminary reduction in \cite[Section 5]{Tru26}. 
Let $\mathfrak b_k$ be the base ideal of $\abs{k\overline{\cL}}$. Take an equivariant log resolution $p_k\colon \cY_k\to\cX$, write $\mathfrak b_k\cO_{\cY_k}=\cO_{\cY_k}(-F_k)$, and define  
 $\cM_k:=p_k^*\cL-\frac{1}{k}F_k$. 
The pair $(\cY_k, \cM_k)$ is a semiample test configuration of $L$. 
The model metric of $(\cY_k, \cM_k)$ is 
\begin{equation}\label{eq:base-ideal-metric}
 \chi_k=f+\frac{1}{k}\log\abs{\mathfrak{b}_k},
 \qquad \log\abs{\mathfrak{b}_k}(v):=-G(v)(\mathfrak b_k).
\end{equation}
See \cite[section 5.4]{BJ22}. 
Multiplication of sections gives
$\mathfrak b_{k_0}^\ell\subset\mathfrak b_{ k_0 \ell}$, so \eqref{eq:base-ideal-metric} implies $\chi_{k_0}\leq \chi_{k_0 \ell}$. Since $\chi_k\leq f$, we have  $\chi_k\leq P(f)=\psi$. 
Let us fix a sufficiently divisible $k_0$ and set $k_j:=k_0 j!$. 
We consequently have
\begin{equation*}
 \chi_{k_0}\leq\chi_{k_j}\leq\psi.
\end{equation*}
In fact $\chi_{k_j}\to\psi$ as $j\to \infty$, uniformly. This is the base-ideal approximation of Boucksom--Jonsson \cite[Lemma 5.22]{BJ22}. For the relative base ideals, it states
\begin{equation*}
 P(f)-f=\lim_{k \to \infty}\frac{1}{k}\log\abs{\mathfrak{b}_k},
\end{equation*}
uniformly when $P(f)$ is continuous. 
This is a non-Archimedean analogue of 
\cite[Theorem 1.5]{Ber09}. 
In particular, we obtain $d_1^{\NA}(\chi_{k_j},\psi)\to0$ and weak convergence of the Monge--Amp\`ere measures.

For these base-ideal approximants, the argument in
\cite[Section 5, Steps 3--4]{Tru26} gives 
\begin{equation}\label{eq:special-fujita}
 \overline{\cM}_k^{n+1}\longrightarrow\vol(\overline{\cL}),\qquad
 \overline{\cM}_k^nK_{\overline{\cY}_k/\overline{\cX}}\longrightarrow0.
\end{equation}
Let us check entropy convergence for the same sequence. Let 
\begin{equation*}
g(v) := G(v)(K_{\overline{\cX}/(X\times \mathbb{P}^1)}+\cX_{0, \red}-\cX_0 ), 
\end{equation*}
where the evaluated divisor is the relative log canonical divisor over $X\times\mP^1$. 
For a divisorial $v$, one has 
$A_X(v)-g(v)=A_{(\overline{\cX}, \cX_{0, \red})}(G(v))$ 
so log canonicity of $(\overline{\cX}, \cX_{0, \red})$ yields $A_X(v)-g(v) \geq 0$. For a component $F$ of $(\cY_k)_0$ and $b_F=\ord_F((\cY_k)_0)$, 
log-smoothness of $(\overline{\cX}, \cX_{0, \red})$ gives 
\begin{equation*}
 0\leq A_X(v_F)-g(v_F)
 =\frac{\ord_F(K_{\overline{\cY}_k/\overline{\cX}}+(\cY_k)_{0,\red}
                         -p_k^*\cX_{0,\red})}{b_F}
 \leq\frac{\ord_FK_{\overline{\cY}_k/\overline{\cX}}}{b_F}.
\end{equation*}
The last inequality uses 
$p_k^*\cX_{0,\red}\geq(\cY_k)_{0,\red}$.
Therefore, integrating by $\MA^{\NA}(\chi_k)$   
\begin{equation*}
 0\leq H^{\NA}(\chi_k)-\int g\MA^{\NA}(\chi_k)
 \leq V^{-1}\overline{\cM}_k^n\cdot K_{\overline{\cY}_k/\overline{\cX}}\longrightarrow0.
\end{equation*}
The limiting measure $\MA^{\NA}(\psi)
=\lim_j \MA^{\NA}(\chi_{k_j})$ is supported at component valuations of $\cX_0$ by \cite[Lemma 4.8]{BJ25}, where $g=A_X$. Continuity of $g$ proves entropy convergence.

We next modify $\chi_{k_j}$ to obtain ample models $\psi_j$. 
On each $\cY_k$ choose a small, $p_k$-exceptional and vertical,  effective divisor $G_k$ such that $\cA_k:=p_k^*\cA_{\cX}-G_k$ is ample. The rational convex combinations
\begin{equation*}
 \cM_{k,s}=s\cM_k+(1-s)\cA_k
\end{equation*}
define an ample test configuration of $(X, L)$. 
For each fixed $k$, their model metrics converge uniformly to $\chi_k$, as $s\to1$.
From the intersection number description, entropy is continuous in $s$. 
Choose rational $s_j$ so that these errors are in $\pm 1/j$. 
Now we define $(\cX_j, \cL_j):=(\cY_{k_j}, \cM_{k_j,s_j})$ and  $\psi_j:=\varphi_{(\cX_j, \cL_j)}$. By the construction we have 
\begin{equation*}
 \chi_{k_0}-1/j\leq\psi_j\leq\psi+1/j.
\end{equation*}
This proves the third assertion of \eqref{eq:model-approximation}. 

As we have already explained, 
associated weak geodesic rays $\ell_{\psi_j}$ 
recovers $(\cX_j, \cL_j)$ by \eqref{eq:ray to NA}. 
The uniform convergence $\abs{\psi -\psi_j}<\e$ is inherited to $\abs{\ell_{\psi,t} -\ell_{\psi_j,t}} <\e t$, 
as we chose the common initial K\"ahler potential. 
It in particular implies $d_p^c(\ell_{\psi}, \ell_{\psi_j}) <\e$ and hence entropy convergence implies 
convergence in $M^{\NA}$. 
Finally, convexity of distance gives
\begin{equation*}
 \abs{E(\ell_{\psi, 1})-E(\ell_{\psi_j, 1})}\leq d_1(\ell_{\psi,1},\ell_{\psi_j,1})
 \leq  d_2(\ell_{\psi,1},\ell_{\psi_j,1})\leq d_2^c(\ell_\psi,\ell_{\psi_j})\longrightarrow0. 
 \qquad
\end{equation*}
Therefore the both terms in the right-hand side of 
\begin{equation*}
 \norm{\ell_{\psi_j} -E^{\rad}(\ell_{\psi_j})}_2^2
 =d_2(\ell_{\psi_j,0},\ell_{\psi_j,1})^2-\bigl(E(\ell_{\psi_j,1})-E(\ell_{\psi_j,0})\bigr)^2.
\end{equation*}
converge and the last assertion follows. 
\end{proof}

\subsection{Approximation of general finite-energy geodesic rays}
\begin{prop}\label{prop:ray-approximation}
Every $\cE^2$-geodesic ray $\ell$ 
emanating from a smooth
potential $u_0\in\cH$,  
with finite radial K-energy, 
admits normal ample test configurations with model metrics $\psi_j$, associated weak geodesic ray $\ell_{\psi_j}$, which satisfy 
\begin{equation*}
 d_2^c(\ell_{\psi_j},\ell)\longrightarrow0,\qquad
 \limsup_{j \to \infty} M^{\NA}(\psi_j)\leq M^{\rad}(\ell).
\end{equation*}
Their $L^2$-norms converge to the centered speed of $\ell$. Namely, 
\begin{equation*}
\norm{(\cX_j, \cL_j)}_2 \longrightarrow \norm{\ell -E^{\rad}(\ell)}_2. 
\end{equation*}
\end{prop}
\begin{proof}
First let $\varphi\in\cE^{1,\NA}$ be bounded, normalized by $\sup_{X^{\NA}}\varphi=0$, and of finite entropy. Choose $M$ with $-M\leq\varphi\leq0$.  Let $\varphi_j$ be a solution to \eqref{eq:retracted-measures} with normalization $\sup_{X^{\NA}}\varphi_j=0$. 
From \cite[Theorem A.1,Proposition 4.12(i)]{BJ22}, 
Monge--Amp\`ere equation \eqref{eq:retracted-measures} implies
\begin{equation*}
 e_p(\varphi_j)
 =\int_{X^{\NA}}(-\varphi_j\circ r_j)^p \MA^{\NA}(\varphi)
 \leq\int_{X^{\NA}}(-\varphi_j)^p\MA^{\NA}(\varphi).
\end{equation*}
In the last inequality we used the fact that retraction maps non-decreases every non-Archimedean plurisubharmonic function: $\varphi_j\leq\varphi_j\circ r_j$. 
We apply Lemma \ref{lem:mixed-moment} with $p=n+2$. 

Suppose $e_p(\varphi_j)\geq1$,  
$M':=\sum_{k=1}^n\binom nk(p)_kM^k$ assures 
\begin{equation*}
 e_p(\varphi_j)\leq M' e_p(\varphi_j)^{\frac{p-1}{p}}. 
\end{equation*}
It implies $e_p(\varphi_j)\leq\max\{1,A^p\}$. Lemma \ref{lem:radial-moment} gives uniformly bounded radial $d_p$ speeds. 
Proposition \ref{prop:retraction} $(1)$ states $d_1^{\NA}$ convergence of $\varphi_j$. 
It combined with the interpolation \eqref{eq:interpolation} now implies $d_2^c$ convergence of their maximal weak geodesic rays. 

Next, for each $j$ we apply Proposition \ref{prop:retraction} $(2)$. 
As $P(f_i\circ r_j)$ uniformly converge to $\varphi_j$, 
The maximal rays with a common initial potential satisfies  $d_2^c(\ell_{P(f_i\circ r_j)}, \ell_{\varphi_j}) \to 0$.  
Apply Proposition \ref{prop:model-approximation} to each envelope $P(f_i \circ r_j)$. 
In each step of approximations we have checked $d_2^c$-convergence of the maximal weak geodesic rays. 
Now diagonal argument yields a subsequence $\psi_k$ which 
converge to $\varphi$ in $d_1^{\NA}$, $d_2^c$, and in entropy. 
Non-entropy term of $M^{\NA}$ is controlled by $d_1^{\NA}$ so convergence in $M^{\NA}$ also 
follows. 

Finally, let $\ell$ be an arbitrary 
$\cE^2$-geodesic ray emanating from $u_0 \in \cH$, with finite radial K-energy. 
Darvas--Lu
\cite[Theorem 1.5(i)]{DL20} established existence of bounded-potential rays $\ell_m$ 
with $\ell_{m, 0}=u_0$, $d_2^c(\ell_m,\ell)\to0$ and 
$M^{\rad}(\ell_m)\to M^{\rad}(\ell)$. 
Moreover, for each $m$, the ray $\ell_m$ is Lipschitz in time. 
That is, there exists a constant 
$C_m$ such that
\begin{equation*}
 u_0-C_mt\leq\ell_{m, t}\leq u_0+C_mt 
\end{equation*}
holds uniformly in $t$. 
Taking the associated non-Archimedean metric \eqref{eq:ray to NA} yields 
\begin{equation*}
 -C_m\leq\varphi_m:=(\ell_m)^{\NA}\leq C_m.
\end{equation*}
For all sufficiently large $m$,
$M^{\rad}(\ell_m)$ is finite. 
The first assertion of Proposition \ref{prop:na-slope} therefore implies that $\ell_m$ is maximal.
Proposition \ref{prop:na-slope} gives
\begin{equation*}
 H^{\NA}(\varphi_m)+R^{\NA}(\varphi_m)
       +\bar S E^{\NA}(\varphi_m)
 \leq M^{\rad}(\ell_m)<\infty.
\end{equation*}
It implies $H^{\NA}(\varphi_m)<\infty$, since 
non-entropy terms are finite for bounded $\varphi_m$. 
Now we apply the first and second step to bounded $\varphi_m-\sup_{X^{\NA}}\varphi_m=:\varphi \in \cE^{1, \NA}$,   take a diagonal sequence, and add back the same constant $\sup_{X^{\NA}}\varphi_m$. 
The associated maximal weak geodesic ray of the obtained sequence $\psi_j:=\psi_{m_j, k_j}$ 
converges to $\ell$ in $d_2^c$ and 
\begin{equation*}
 \limsup_{j \to \infty} M^{\NA}(\psi_j)\leq M^{\rad}(\ell). 
\end{equation*}
Since $E(\ell_{1})-E(\ell_{0})=E^{\rad}(\ell)$ for all $\cE^2$-geodesic rays and 
\begin{equation*}
\abs{E^{\rad}(\ell) -E^{\rad}(\ell_{\psi_j})} 
\leq d_1^c(\ell, \ell_{\psi_j}) \leq d_2^c(\ell, \ell_{\psi_j}) \longrightarrow 0, 
\end{equation*}
The $L^2$-norm 
\begin{equation*}
\norm{(\cX_j, \cL_j)}_2^2=
 \norm{\ell_{\psi_j} -E^{\rad}(\ell_{\psi_j})}_2^2
 =d_2(\ell_{\psi_j,0},\ell_{\psi_j,1})^2-\bigl(E(\ell_{\psi_j,1})-E(\ell_{\psi_j,0})\bigr)^2 
\end{equation*}
converges to $ \norm{\ell -E^{\rad}(\ell)}_2^2$. 
\end{proof}

\section{Completion of the proof and related results in the Fano case}\label{sec:conclusion}

\subsection{Competion of the proof of Theorem \ref{thm:main}}

\begin{proof}[Proof of Theorem \ref{thm:main}]
At the end of section ~\ref{sec:normalization} we proved 
\begin{equation*}
\inf_{\cH}\Cal \geq \inf_{\cE^2}\abs{\partial M} \geq 
\sup \frac{-M^{\NA}(\cX, \cL)}{\norm{(\cX, \cL)}_2}.   
\end{equation*}
Theorem \ref{prop:smooth-recovery} shows the first inequality is equality. It remains to prove that the last inequality is equality. 
Since $(X,L)$ is K-unstable, $\alpha:=\inf_{\cE^2}\abs{\partial M} >0$. Choose a maximizing ray in \eqref{eq:xia}. 
Subtract the slope of affine $E$. Slope of K-energy is unchanged and squared speed decreases by the square of the subtracted scalar. Speed of the subtracted ray is nonzero, otherwise the ray is constant and has zero radial K-energy. The minus ratio of radial K-energy with speed is still at least $\alpha$, and \eqref{eq:xia} bounds it above by $\alpha$; hence it still attains $\alpha$.

By Proposition \ref{prop:ray-approximation}, there are normal ample
test configurations with centered norms $N_j\to N:=d_2(\ell_0, \ell_1)>0$
and $\limsup_j M^{\NA}(\psi_j)\leq M^{\rad}(\ell)=-\alpha N$.
For every $\delta>0$, eventually
$M^{\NA}(\psi_j)\leq-\alpha N+\delta$, whence
\begin{equation*}
 \liminf_j\frac{-M^{\NA}(\psi_j)}{N_j}\geq\alpha.
\end{equation*}

Since $N_j\to N>0$, we have $N_j>0$ for all sufficiently
large $j$.
The already established inequality gives
$\frac{-M^{\NA}(\psi_j)}{N_j}\leq\alpha$. 
Thus the quotients converge to $\alpha$,
which proves \eqref{eq:main}.

Moreover, the above argument gives
$M^{\NA}(\psi_j)\geq-\alpha N_j$.
Together with $N_j\to N$ and
$\limsup_jM^{\NA}(\psi_j)\leq-\alpha N$,
this proves
$M^{\NA}(\psi_j)\to-\alpha N=M^{\rad}(\ell)$. 
\end{proof}

\subsection{Related results in the Fano case}

We close this article by touching upon some related results in the Fano case. 
For the anti-canonical polarization $L=-K_X$, the cscK condition is characterized by the vanishing of normalized Ricci potential 
$\Ric(\omega_\varphi)-\omega_\varphi=dd^c\rho_\varphi$, $V^{-1}\int_Xe^{\rho_\varphi}\omega_\varphi^n=1$. 
The moment--weight equality for the Ricci analogue of Calabi functional 
\begin{equation*}
\left[\frac{1}{V}\int_X(e^{\rho_\varphi}-1)^2\omega_\varphi^n\right]^{\frac{1}{2}}
\end{equation*}
and the corresponding non-Archimedean Ding energy was proved in \cite[Theorem A]{His19}. 
It is interesting to compare the infima of 
these two Calabi-type functionals. 

In \cite[Section 4]{HI26}, the authors explain how
Donaldson's conjecture implies a Miyaoka--Yau inequality
incorporating the delta invariant.   
Combined with Theorem \ref{thm:main}, we obtain the following. 
The same inequality is established also for klt Fano
varieties \cite[Theorem 1.1]{HI26}, by a different approach. 

\begin{cor}\label{cor:my}
Every smooth Fano manifold of dimension $n\geq2$ satisfies
\begin{equation*}
 \left(2(n+1)c_2(X)-nc_1(X)^2\right)c_1(X)^{n-2}
 \geq-n\left(1-\min\{1,\delta(X)\}\right)^2c_1(X)^n.
\end{equation*}
\end{cor}

As it was explained in \cite{HI26}, 
the present approach 
clarifies a differential geometric meaning of the equality case. 

\bibliographystyle{amsplain}
\bibliography{calabi-mabuchi}
\end{document}